# Design and Operation of Renewable Energy Microgrids under uncertainty towards Green Deal and Minimum Carbon Emissions

Su Meyra Tatar[a], Erdal Aydin[a,b,c*]

[a]Department of Chemical and Biological Engineering, Koç University, Istanbul 34457, Türkiye
[b]Department of Chemical Engineering, Bogazici University, Bebek, Istanbul 34342, Türkiye
[c]Koç University TUPRAS Energy Center (KUTEM), Koç University, Istanbul, 34450, Türkiye

*eaydin@ku.edu.tr

**Abstract**

The regulations regarding the Paris Agreement are planned to be adapted soon to keep the global temperature rise within 2 $^0$C. Additionally, integrating renewable energy-based equipment and adopting new ways of producing energy resources, for example *Power to Gas* technology, becomes essential because of the current environmental and political concerns. Moreover, it is vital to supply the growing energy demand with the increasing population. Uncertainty must be considered in the transition phase since parameters regarding the electricity demand, carbon tax policies, and intermittency of renewable energy-based equipment have intermittent nature.

A multi-period two-stage stochastic MILP model is proposed in this work where the wind speed, solar irradiance, temperature, power demand, carbon emission trading (CET) price, and $CO_2$ emission limit are considered uncertain parameters. This model finds one single optimal design for the energy grid while considering several scenarios regarding uncertainties simultaneously. Three stochastic case studies with scenarios including different combinations of the aforementioned uncertain parameters are investigated. Results show that more renewable energy-based equipment with higher rated power values is chosen as the sanctions get stricter. In addition, the optimality of PtG technology is also investigated for a specific location. Implementing the $CO_2$ emission limit as an uncertain parameter instead of including CET price as an uncertain parameter results in lower annual $CO_2$ emission rates and higher net present cost values.

## 1. Introduction

### 1.1. Motivation

Today, we are in a situation where we extremely observe the effects of global warming [1–4]. The main concern of today's world is to prevent global warming while keeping the global temperature increase within $2^0C$. The sanctions that come with Paris Agreement are essential and binding for energy producers where the ultimate aim of countries is to reach net zero emission in the near future [5]. Another concern of energy producers is to meet the enhanced energy demand with the increasing population. These concerns drive the energy producers to supply this demand and decarbonize the energy industry at the same time [6,7]. Due to the aforementioned challenges in the energy supply industry, the integration of renewable energy resources into the energy grids is essential.

To decide on the design of a new generation energy grid from today that is planned to operate for many years, it is required to foresee several parameters from this day, which is an extremely challenging task. Hence, instead of formulating deterministic models that only hold for one certain case, stochastic programming helps to integrate the uncertain nature of parameters into the decision-making model. Regarding the sanctions that come with the Paris Agreement, uncertain parameters corresponding to the renewable energy intermittency and carbon emission policies are the most crucial ones.

### 1.2. Background

Ref. [8] gives insights into the uncertain parameters having different time scales in energy systems. While energy demand and price, and energy produced by renewables have short and medium-term time scales, uncertainties regarding technological developments and policies in the energy sector,

It is crucial to represent long and short-term time scale uncertainties when designing and planning an energy hub. MILP models are used widely to formulate multi-input multi-output energy systems, and several studies in the literature consider some of the above uncertainties. For example, in [9], the uncertain parameters, which are electricity demands, price of electricity, and natural gas, vary according to the four seasons. In another study [10], electricity demand and supply are uncertain parameters. In [11], the design and operation of a renewable energy system are investigated via formulating the problem as a two-stage stochastic MILP. An energy system with wind turbines, photovoltaics (PV), and energy storage system investigates the five different scenario generation methods in case studies. In [12], the optimal operation of the seaports is studied while formulating the energy management system as MILP together with considering the uncertain parameter of renewable energy generation. The model given in Ref. [13] is modeled as multi-stage stochastic MILP, which considers the hydrogen demand as an uncertain parameter, and the design of the energy system is investigated with two case studies; different ways of producing hydrogen. In [13], the uncertain parameters are determined according to their higher effect on the objective function. Hydrogen demand is chosen as an uncertain parameter in that study due to its higher impact. Although carbon tax is in second place with the second highest effect, it is not included in the study.

Although these studies consider some of the aforementioned uncertain parameters, they do not consider the $CO_2$ limitations and focus on the design and operation of the energy grid simultaneously. In addition, none of them study renewable energy integration and PtG technology altogether, although the formulated problems are stochastic. Accordingly, to obtain more realistic results considering the current trends related to carbon emission taxing and cap and trade dictated by the Green Deal, it is vital to integrate uncertainties regarding the carbon pricing policies into a model for optimal design and operations of energy systems.

In reducing carbon emissions, it is important to price carbon properly while accounting for environmental concerns. The two most common terms regarding carbon price are carbon allowances and carbon tax [14,15]. Carbon allowances, 'cap and trade' in other words, refer to a unit of carbon dioxide that can be emitted. This amount might be directly given or purchased depending on national government sanctions. In 'cap and trade,' while emitting less, 'cap' can be traded with a price where this system drives the firms/governments to emit less. On the other hand, in the carbon tax approach, a tax is applied to firms for emitting excess amounts of $CO_2$. The main difference between these two approaches is as follows; in the carbon allowance approach, allowable $CO_2$ emission is limited while the carbon trading price may change; in the carbon tax approach, the tax is constant while the emitted amount of $CO_2$ may change [14,16]. In this work, the carbon allowance approach is considered, and carbon emission is limited per the Paris Agreement, where carbon emission trading (CET) price and $CO_2$ emission limit are considered uncertain parameters.

There are several studies concerning emission reduction without considering the uncertain nature of parameters regarding carbon emission. In [6], a renewable energy system includes wind turbines and photovoltaic panels to meet a certain location's heat and electricity demand. This study aims to make greenhouse gas emissions nearly zero for a city in Italy by 2030. In [7], carbon neutrality is aimed to reach 2050 in Switzerland with biogenic carbon capture, where the model is formulated as MILP.

There are also some papers in the literature that consider the $CO_2$ limitations but containing only a few of the uncertain parameters mentioned previously. In another study [10], power demand and supply are uncertain parameters, and this work investigates the optimal design of an energy system while formulating the problem as a two-stage stochastic MILP without considering the carbon emission constraints. Conversely, a previous work [17] considers the emission limitations, but it investigates the optimal design only deterministically. Ref. [18] modeled the multi-network problem as a two-stage stochastic MILP with the uncertainties of the product price, demand, and raw material price, where this study aims to minimize carbon emissions.

### 1.3. Proposed Work

In this work, a two-stage stochastic MILP model is formulated to find the optimal design and the sequence of equipment in the energy hub. The stochastic MILP model is formulated to meet the one-third power and heat demand of a certain location in Türkiye. The hub superstructure is integrated with renewable energy-based equipment and PtG technology, focusing on the importance of hydrogen, which considers several vital uncertain parameters simultaneously. Six uncertain parameters are included, and different combinations of them are investigated via several case studies. The uncertain parameters are power demand and the ones related to the intermittency of renewable energy-based equipment and $CO_2$ emission constraints. Wind speed, solar irradiance, and temperature data are the ones regarding the intermittency of renewables. The ones regarding the emissions are carbon emission trading (CET) price and $CO_2$ emission limit. Including several uncertainties usually increases the computational complexity drastically [13]. To keep the computational complexity in an acceptable region, sixteen candidate pieces of equipment are included in the equipment pool for the first case study. Other case studies include an additional electrolyzer and methanation reactor representing the PtG technology and the impact of green hydrogen [19]. A multi-period two-stage stochastic MILP model is formulated to represent a multi-input, multi-output energy grid which is integrated with renewable energy-based equipment and PtG technology under uncertainty.

Three case studies that include different combinations of these uncertain parameters are studied to observe the effect of these parameters. Additionally, $CO_2$ emission is limited according to Paris Agreement, and in consideration of the sanctions regarding the Paris Agreement, CET price and $CO_2$ emission limit are introduced as uncertain parameters in different case studies to investigate their effect of them on the optimal decision of installation and operation of the equipment in energy grid for the lifetime of the project.

The organization of the paper is as follows: Section 2 gives insight into stochastic programming. Section 3 includes the mathematical formulation of the model and introduces the uncertain parameters in this study. In Section 4, results and discussion of three case studies are given. Finally, Section 5 is the conclusion of this work.

## 2. Two-stage Stochastic Programming

### 2.1. General Formulation

The general two-stage stochastic programming formulation is given below. Here, *w* represents the scenarios, and *E* represents the mathematical expectation operator for every scenario, *w* [21].

$$\min_{x} c^T x + E_w Q(x, w)$$

$$\text{s.t.} \quad Ax = b,$$

$$x \geq 0,$$

$$Q(x, w) \equiv \min_{y} (d_w)^T y$$

$$\text{s.t.} \quad T_w x + W_w y = h_w,$$

$$y \geq 0, \quad \forall w \in \Omega,$$

Decision variables are divided into the first and second-stage decision variables in a two-stage mixed-integer linear programming formulation. Here, while the x vector represents the first stage decision variables, the y vector represents the second stage decision variables. The second stage decisions are the ones that have been made after the first stage decisions are realized [20]. The first stage decision variables are related to the design stage; and are called 'here and now' decision variables, whereas the second stage variables are called 'wait and see' recourse variables, and they are related to the operational stage [10].

### 2.2. Deterministic Equivalent

The deterministic equivalent of the general formulation of the two-stage stochastic programs is given below, which is written for finite number uncertainty realizations where *p* represents the probability of occurrence for every *w*. It should be noted that *x* and *y* decision variables can be discrete and continuous.

$$\min_{x, y_w \, \forall w} c^T x + \sum_w p(w)(d_w)^T y_w$$

$$\text{s.t.} \quad Ax = b,$$

$$T_w x + W_w y_w = h_w, \forall w,$$

$$x \geq 0, y_w \geq 0.$$

In this formulation, the capacity of storage units, the rated power of the equipment, and the binary decision variables to install equipment are the first stage decision variables; in other words, design decision variables.

## 3. Proposed Stochastic MILP Model

### 3.1. Scenario Generation & Uncertain Parameters

In two-stage stochastic optimization, the decision-making model results in a single optimal solution for the design phase considering all possible scenarios. On the other hand, in deterministic optimization, optimal design and operation are found according to one single scenario, which is one of the main advantages of stochastic programming [22].

This work considers three base scenarios namely, Likely, Midlikely and Unlikely cases with different demand values, wind speed, temperature, and solar irradiance data. All scenarios utilized in the stochastic formulations are derived from these base scenarios. Please note that the base scenarios used in this work are generated from [23]. The demand data were taken from the official electricity distributor of Yalova, UEDAS. In addition, all weather-related data, e.g., wind speed, solar irradiation, ambient temperature, etc., are obtained from actual data using an API. The wind speed, solar irradiance, and electricity demand values for Likely, Midlikely, and Unlikely cases are given in Figs. 1-3, respectively. The six uncertain parameters that are included in this study are given in Table 1, which represents all possible values of the uncertain parameters. It should be noted that the scenarios in case studies are produced from the combination of the values of these uncertain parameters.

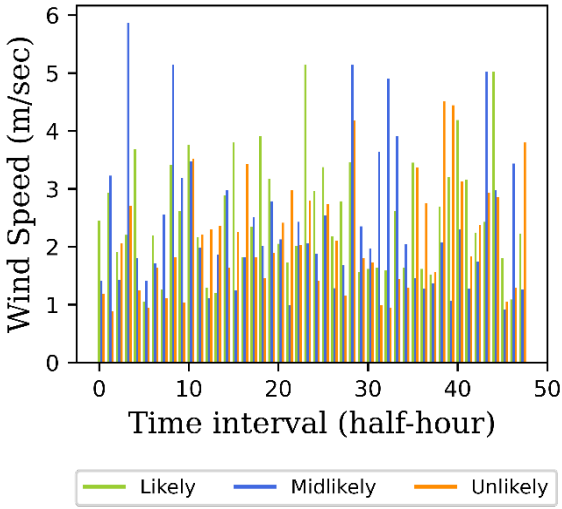

**Figure 1.** Wind speed data for likely, mid-likely and unlikely scenarios.

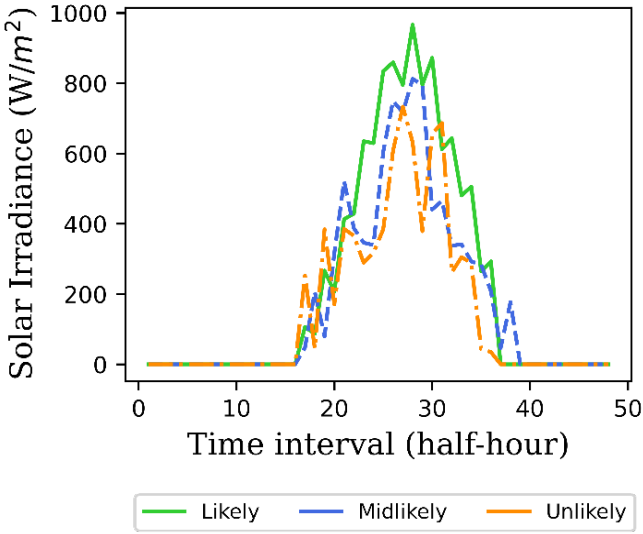

**Figure 2.** Solar irradiance data for likely, mid-likely and unlikely scenarios.

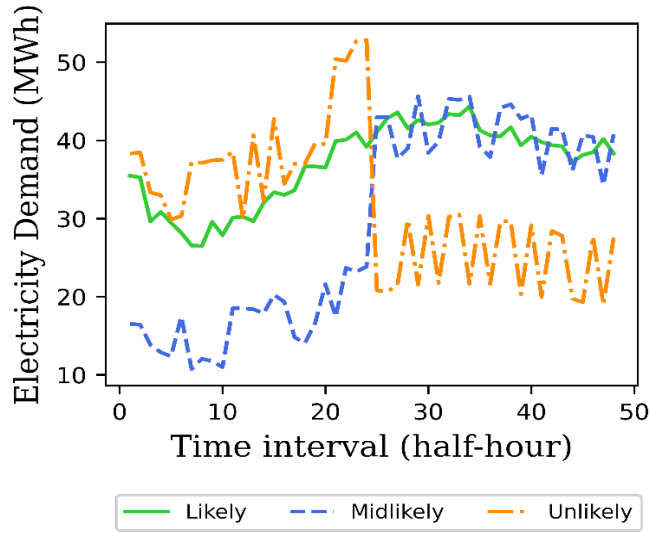

**Figure 3.** Electricity demand data for likely, mid-likely and unlikely scenarios.

**Table 1.** The possible values of all uncertain parameters where the scenarios in case studies are produced from.

|  | Operating Power Coefficient of Renewables as the fraction of the Rated Power | | | | |
| --- | --- | --- | --- | --- | --- |
| **Electricity Demand** | **Wind Speed** | **Solar Irradiance** | **Temperature** | **Carbon Emission Trading Price** | **$CO_2$ Emission Limit** |
| Likely Case | Likely Case | Likely Case | Likely Case | Increasing | Decreasing |
| Midlikely Case | Midlikely Case | Midlikely Case | Midlikely Case | Constant | Constant |
| Unlikely Cases | Unlikely Case | Unlikely Case | Unlikely Case | Slightly Increasing | Slightly Decreasing |

Figs. 4 and 5 represent the operating power coefficient of wind turbines and solar cells, respectively, as the fraction of the rated power for likely, mid-likely, and unlikely cases.

In Fig. 4, it is observed that the wind availability of scenario $w_2$ is the highest, while scenario $w_3$ has the lowest wind availability among other scenarios.

While wind availability depends on wind speed, PV availability depends on temperature and solar irradiance data, as given by Eqs. 24 and 25. The PV availability of likely case is the highest among other scenarios. Although PV availabilities of mid-likely and unlikely cases are seen very close to each other, the PV availability of mid-likely case is greater than the unlikely case. It should be noted the wind and solar availability comparison of scenarios can be observed from the power balance figures of the case studies when wind and solar cells are chosen by the optimal decision-making model.

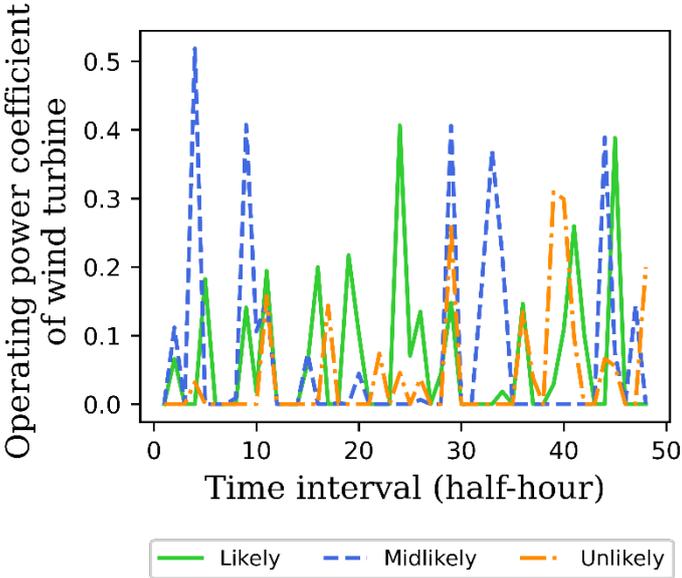

**Figure 4.** Operating power coefficient of wind turbine as the fraction of the rated power for three cases.

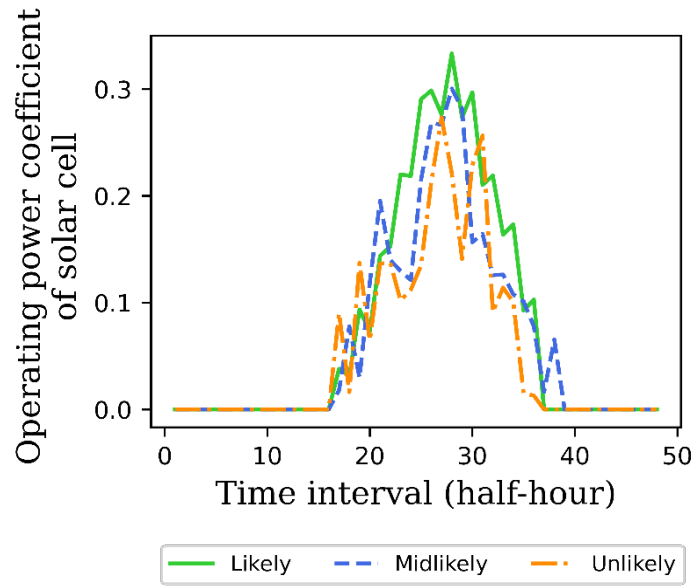

**Figure 5.** Operating power coefficient of PV as the fraction of the rated power for three cases.

It is aforementioned that limiting and controlling $CO_2$ emission is essential to prevent future environmental damages. Hence, it is important to price the carbon emission and determine the $CO_2$ limit properly [16]. Since pricing $CO_2$ emission to a certain value and assigning a certain $CO_2$ limit to the model does not represent the needs for a 20-year project, carbon emission trading (CET) price and $CO_2$ emission are taken as uncertain parameters and investigated in Case studies 2 and 3, respectively.

### 3.2. Mathematical Formulation of the Model

The mathematical formulation of the two-stage stochastic MILP model is represented in the supplementary materials section due to space limitations. For detailed explanations, the reader is referred to the appendix section.

## 4. Case Studies

The two-stage stochastic MILP decision-making model is coded in GAMS language, and the CPLEX solver is used. Twenty-year of project lifetime is considered and the lifetimes of equipment are assumed to be consistent. Demand and carbon tax profiles change in each five years as described in the following sections. It is also possible to add seasonality to each year representation, which on the other hand increases the overall complexity of the two-stage stochastic model and would probably require special decomposition techniques.

Case study 1 includes 16 candidate pieces of equipment. Case studies 2 and 3 include an additional compact *Power to Gas* unit, which includes an electrolyzer and a methanation reactor, in which the optimality of whether to install this unit is investigated. The superstructure of the model is given in Fig. 6. For this study, although electricity purchase from the national grid is allowed to some limit, heat purchase is not allowed. It is also assumed that the operational hours of the grid are not limited; it operates 24 hours a day. The differences between case studies are summarized in Table 2, in which detailed explanations are included in case studies.

**Table 2.** Summary of differences of case studies.

|  | Case 1 | Case 2 | Case 3.a | Case 3.b | Case 3.c |
|---|---|---|---|---|---|
| **Increase percentages in CET price** | - | 20%<br>10%<br>0% | - | - | - |
| **Decrease percentages in $CO_2$ emission limit** | - | - | 20%<br>10%<br>0% | 50%<br>20%<br>0% | 50%<br>20%<br>0% |
| **Terms in Objective Function** | Cost Stream | Cost Stream- ('Cap and Trade' Income +' Excess SNG' Income) | Cost Stream- (' Excess SNG' Income) | Cost Stream- (' Excess SNG' Income) | Cost Stream- (' Excess SNG' Income) |

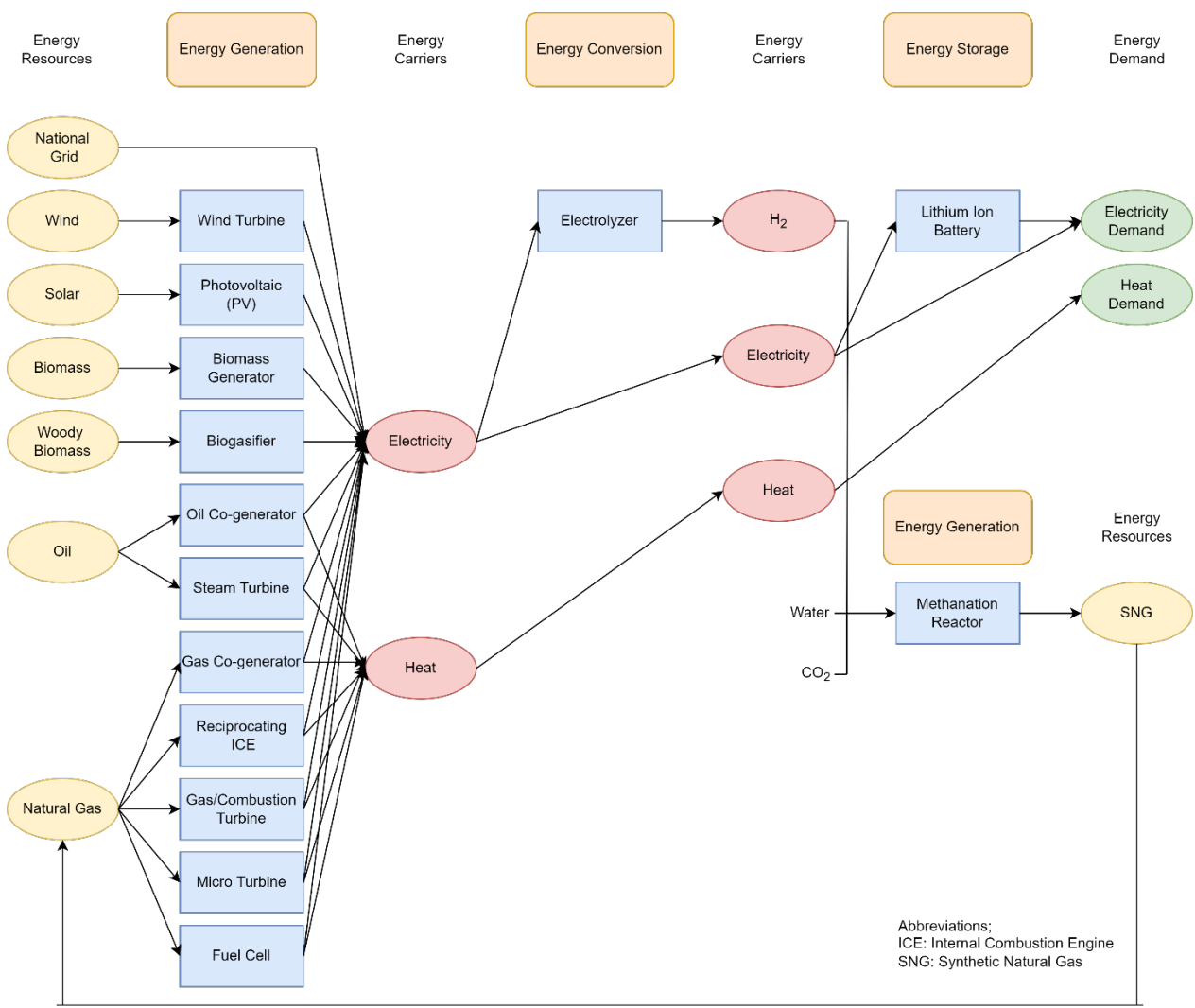

**Figure 6.** Superstructure of the stochastic MILP model.

## 4.1. Case 1

Case 1 intends to investigate the optimal equipment and scheduling sequence choice for a specific location while considering uncertainties regarding base scenarios. It is also aimed to observe the gains in stochastic formulation over the deterministic one while calculating the value of the stochastic solution (VSS). Additionally, carbon emission trading (CET) price and $CO_2$ emission limit is taken as uncertain parameters over the 20-year period and introduced as additional uncertain parameters in Case studies 2 and 3, respectively.

The scenarios in Case 1 are represented as $w_1$, $w_2$, and $w_3$, respectively. Please note that every scenario has an equal probability of 1/3 [23]. The scenarios and the corresponding values of uncertain parameters in Case 1 are given in Table 3, in which the corresponding data are given Figs. 1-3 in Section 3. It should be noted that the Case 1 includes uncertain parameters of electricity demand and operating power coefficient of renewables as the fraction of the rated power. Other uncertain parameters given in Table 1 are not included in Case 1. Three different scenarios for Case 1 are produced from possible values of four of the uncertain parameters out of six given in Table 1.

Table 3. Uncertainties for Case 1.

| Scenarios | Electricity Demand | Operating Power Coefficient of Renewables as the Fraction of the Rated Power | | | Probability |
|---|---|---|---|---|---|
| | | Wind Speed | Solar Irradiance | Temperature | |
| $w_1$ | Likely Case | Likely Case | Likely Case | Likely Case | 1/3 |
| $w_2$ | Midlikely Case | Midlikely Case | Midlikely Case | Midlikely Case | 1/3 |
| $w_3$ | Unlikely Case | Unlikely Case | Unlikely Case | Unlikely Case | 1/3 |

The optimal result of the two-stage stochastic decision-making model for this case study is found as 11.1 billion TRY over 20 years of plant life when a 12% discount rate is considered. The chosen equipment, which is first-stage variables, with corresponding rated power values and the problem size are shown in Tables 4 and 5, respectively.

Table 4. Chosen Equipment for Case 1.

| Equipment | Rated Power (kW) | Capacity (kWh) |
|---|---|---|
| PV3 | 191886.7 | |
| OilCogen1 | 9669.3 | |
| Battery | 24048.5 | 100000 |
| RecipEngine3 | 9341 | |
| FuelCell1 | 1400 | |
| BioGasif1 | 6600 | |

Table 5. Size and CPU time for Case 1.

| | Stochastic Solution |
|---|---|
| CPU Time (sec) | 12431.45 |
| Continuous Variables | 377,403 |
| Discrete Variables | 28,816 |
| Number of Constraints | 643,736 |

The $CO_2$ emission limit is a vital constraint that should be satisfied. Under current green deal related regulations, which are expected to get stricter every day, renewable energy-based equipment is more likely to be chosen by the decision-making model due to their low $CO_2$ emission rates. Especially wind turbines and solar cells may be preferable because they do not require raw material purchase under available wind and sunlight. However, wind turbines are not optimal for the investigated location due to the wind profile. Although solar cells are expensive equipment, they are chosen by the model due to their low emission rates and free of charge raw material. It should be noted from Table 3 that scenario $w_1$ includes the likely case's wind speed, solar irradiance, temperature, and demand data. In the same manner, while scenario $w_2$ has the mid-likely case uncertain parameters, scenario $w_3$ has the unlikely scenario data.

The second-stage recourse variables for each scenario are reported separately. Electricity demand-supply balances are given in Figs. 7-9 for three scenarios. One of the primary factors for the corresponding equipment selection is the low $CO_2$ emission rates. Additionally, temperature and solar irradiance are some factors for PV units to be installed. PV availability comparison among scenarios is given in Fig. 5 and can also be observed in the power demand-supply balance figures, which are given in Figs. 7-9.

For all scenarios, most of the power demand throughout the day is supplied with PV units. Since the PV availability is highest in the likely case, the amount of power supplied by the PV unit is most significant in scenario $w_1$, among other scenarios. Scenarios $w_2$ and $w_3$ follow it since PV availability in the mid-likely case is higher than in the unlikely case. The second and the third highest power supply is realized by a reciprocating combustion engine and oil co-generator, respectively, for scenario $w_1$. For scenario $w_2$, the second and third highest power supply is realized by bio-gasifier and reciprocating combustion engine; for scenario $w_3$, the second and third highest power supply is made by a reciprocating combustion engine and bio-gasifier, respectively. Some of the power demand is supplied from outside of the grid. The power is purchased up to its maximum limit of 5 MWh for most time intervals for all scenarios. Most importantly, the critical advantage of the proposed approach is that even though the recourse variables might differ from scenario to scenario, the stochastic MILP model computes only a single best design for all scenarios.

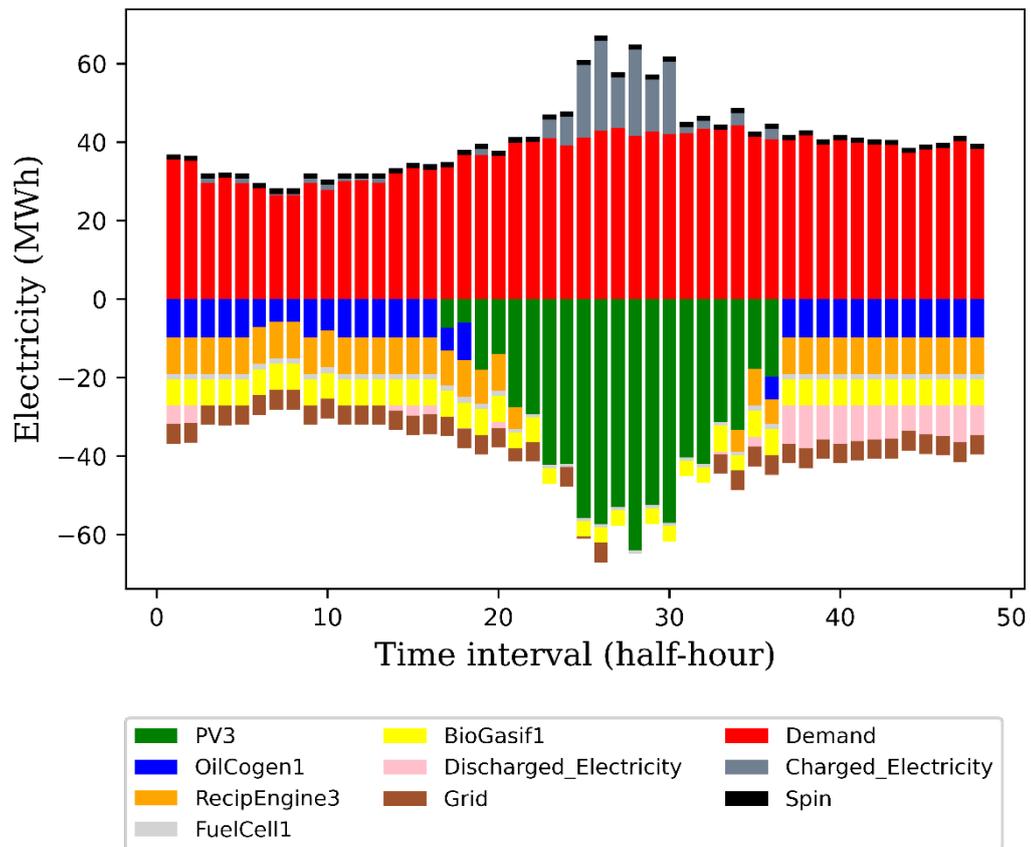

**Figure 7.** Optimal stochastic recourse results for 1st year for scenario $w_1$ in Case 1.

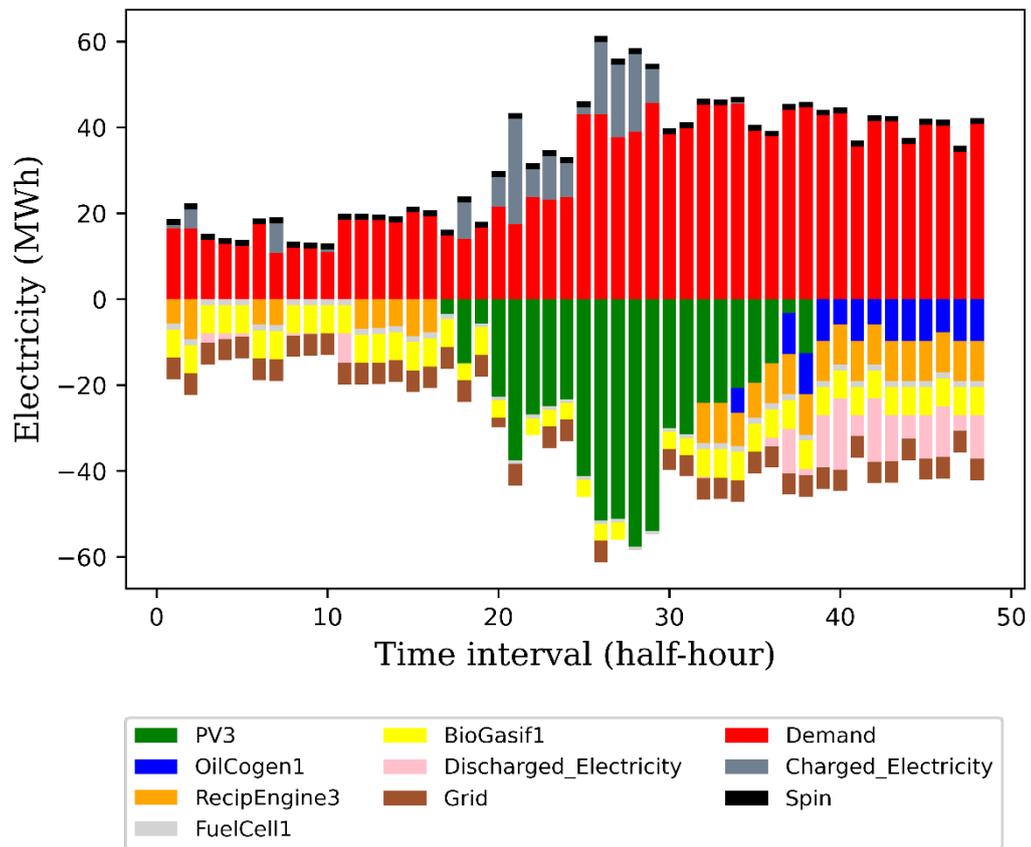

**Figure 8.** Optimal stochastic recourse results for 1$^{st}$ year for scenario $w_2$ in Case 1.

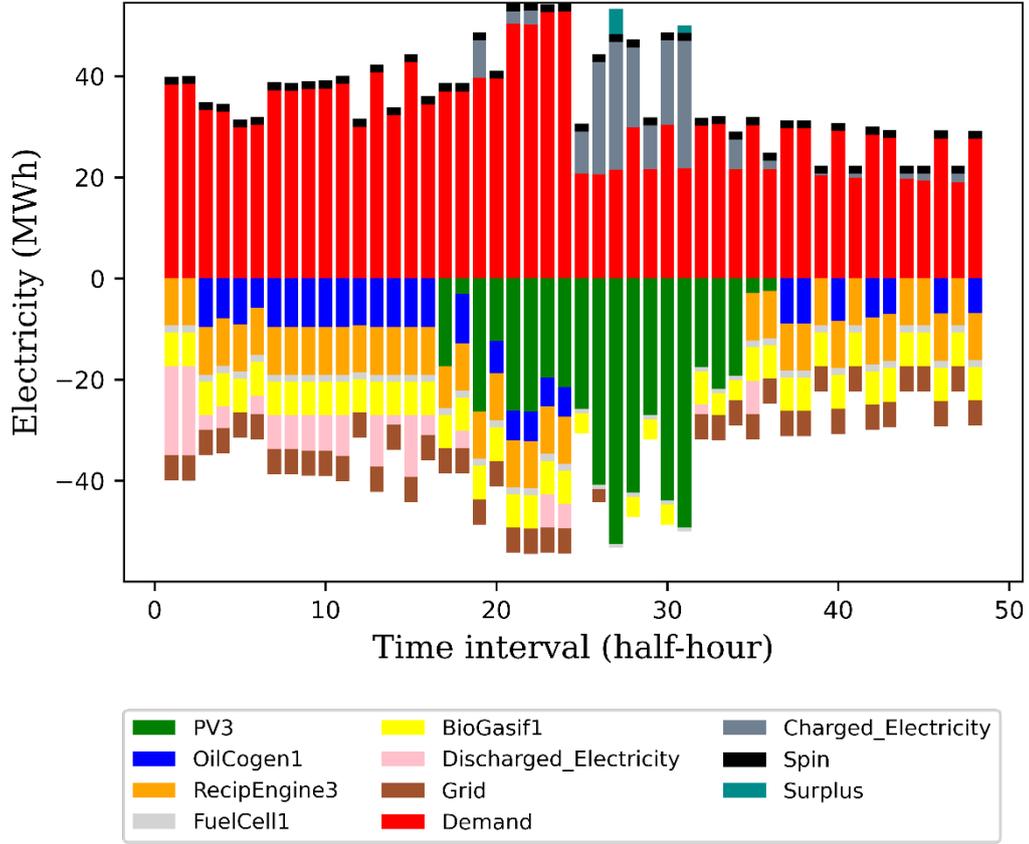

**Figure 9.** Optimal stochastic recourse results for 1st year for scenario $w_3$ in Case 1.

VSS (Value of stochastic solution) index is used to observe the gains in stochastic formulations over the deterministic ones. It typically represents *"how much we can gain when the uncertainties are considered explicitly"* [25]. The VSS index is calculated with Eq. 1, where SS stands for the stochastic solution:

$$\text{VSS} = \text{EVS} - \text{SS} \tag{1}$$

In this formulation, the capacity of storage units, the rated power of the equipment, and the binary decision variables to install equipment are the first stage decision variables; in other words, design decision variables. The cost obtained from SS is 11.1 billion TRY when a 12% discount rate is considered. To calculate the VSS index, the EVS term is also calculated where EVS stands for 'expected value solution' [25]. In calculating the EVS term, first, the deterministic formulation is solved with average values of scenarios' uncertain parameters. Then, the obtained first stage decision variables are fixed. With these fixed design parameters, the deterministic model is solved for recourse variables with every scenario separately. Then, the expected value solution is obtained with the summation of the results obtained from the deterministic model with fixed variables and with corresponding probabilities multiplied. Once EVS is calculated, the VSS index is calculated with Eq. 1.

To illustrate the EVS index calculation, the deterministic model is solved for scenario $w_1$ using the likely case data (likely case's uncertain parameters) with fixed first stage variables, and a cost term is obtained. The decision variable set consists only of recourse variables in this setting. This is repeated for other scenarios as well, and three different cost terms are obtained with fixed design variables. In fact, infeasibilities are obtained at this stage of the calculations, implying that the VSS is $+\infty$ practically. Even obtaining a $+\infty$ VSS value means that the results when considering the uncertainties in the stochastic formulation are way more reliable than the deterministic results. The reason behind the infeasibilities is that the design parameters found for the average demand value are infeasible for scenario $w_1$ because the likely case has a greater demand value than the average demand value of three cases. Another would be for scenario $w_2$, which has a lower demand value than the average demand value of three cases. This case will not be feasible unless a certain amount of excess electricity production is allowed. Therefore, it can be concluded that the two-stage stochastic formulation must be preferred over the deterministic one.

On the other hand, it is still vital to analyze the difference between the deterministic and stochastic versions in a feasible case. Thus, the original problem is relaxed due to the infeasibility. To deal with the infeasibilities, some of the constraints are relaxed for the purpose of obtaining a viable VSS index not equal to infinity. When the constraints regarding the electricity input and the surplus electricity production are relaxed, a 5.14% VSS index is obtained. It may be considered a low percentage at first glance. On the other hand, please note that a substantial amount of manual labor is required to obtain a feasible solution. Furthermore, considering that the net present cost is on the scale of billion liras, a significant amount of savings can be done by stochastic programming based decision-making.

Fig. 10 shows that the demand difference between the three cases and average demand is greatest for the mid-likely case until the 25[th] time interval. After that, the difference is highest for the unlikely case. The time intervals of 13 and 15 are the first two peaks in the demand figure of the unlikely case, as shown in Fig. 3 in Section 3. The optimal solution found with fixed variables (using average uncertain parameters) does not meet the higher demand of the unlikely case. The same applies to the time intervals of 26 and 27. In these time intervals, since the first stage variables are computed for the average demand values, which are higher than the unlikely case's demand, no or poor relaxation of the surplus constraint results in the infeasibility of the power supply equation.

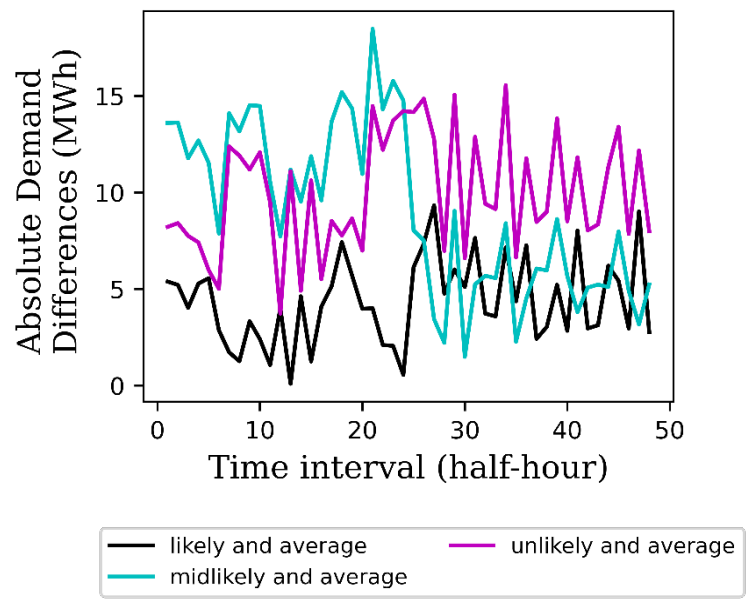

**Figure 10.** Absolute difference values between demand values of scenarios and average demand.

## 4.2. Case 2

In this case study, electricity demand and carbon emission trading (CET) price are considered uncertain parameters in addition to solar irradiance, wind speed, and ambient temperature. Table 6 shows the scenarios of Case 2 containing five uncertain parameters. These scenarios in Case 2 are created from the combination of previous base scenario data and the CET price uncertainty. It should be noted that an equal probability of 1/9 is assigned to every scenario. The sixth uncertain parameter given in Table 1, $CO_2$ emission limit, is not considered an uncertain parameter in this case study. It is taken as constant with a value of $280 gCO_2/kWh$, which is taken as the base value for the $CO_2$ emission limit for this study regarding the Paris agreement.

Table 6. Uncertainties for Case 2.

| Scenarios | Electricity Demand | Operating Power Coefficient of Renewables as the Fraction of the Rated Power | | | CET Price | Probability |
| --- | --- | --- | --- | --- | --- | --- |
| | | Wind Speed | Solar Irradiance | Temperature | | |
| $w_1$ | Likely Case | Likely Case | Likely Case | Likely Case | Increasing | 1/9 |
| $w_2$ | Midlikely Case | Midlikely Case | Midlikely Case | Midlikely Case | Increasing | 1/9 |
| $w_3$ | Unlikely Cases | Unlikely Case | Unlikely Case | Unlikely Case | Increasing | 1/9 |
| $w_4$ | Likely Case | Likely Case | Likely Case | Likely Case | Constant | 1/9 |
| $w_5$ | Midlikely Case | Midlikely Case | Midlikely Case | Midlikely Case | Constant | 1/9 |
| $w_6$ | Unlikely Cases | Unlikely Case | Unlikely Case | Unlikely Case | Constant | 1/9 |
| $w_7$ | Likely Case | Likely Case | Likely Case | Likely Case | Slightly increasing | 1/9 |
| $w_8$ | Midlikely Case | Midlikely Case | Midlikely Case | Midlikely Case | Slightly increasing | 1/9 |
| $w_9$ | Unlikely Cases | Unlikely Case | Unlikely Case | Unlikely Case | Slightly increasing | 1/9 |

The possible values for CET price are given in Table 1 in Section 3. Considering the regulations regarding the Paris Agreement, three possible cases are considered for the value of CET price in the future. It is guessed that the CET price would be constant, slightly increasing, or increasing during the project's lifetime, where decrease in CET price is not even a possibility for the future considering the global warming issue. After choosing a base CET price, the slightly increasing and increasing cases are planned by implementing a 10% and 20% increase in CET price every five years. These increase percentages are chosen to reach and be coherent with the CET price values in the literature that are forecasted for the future, especially for the duration of this project.

It is predicted that the carbon emission trading (CET) price would be constant for three of the scenarios during the lifetime of the proposed work. For the other three scenarios, it would increase slightly so that there will be a 10% increase every five years. Finally, the remaining three scenarios foresee that it will increase by 20% every five years. The predicted carbon emission trading price values are shown in Fig. 11.

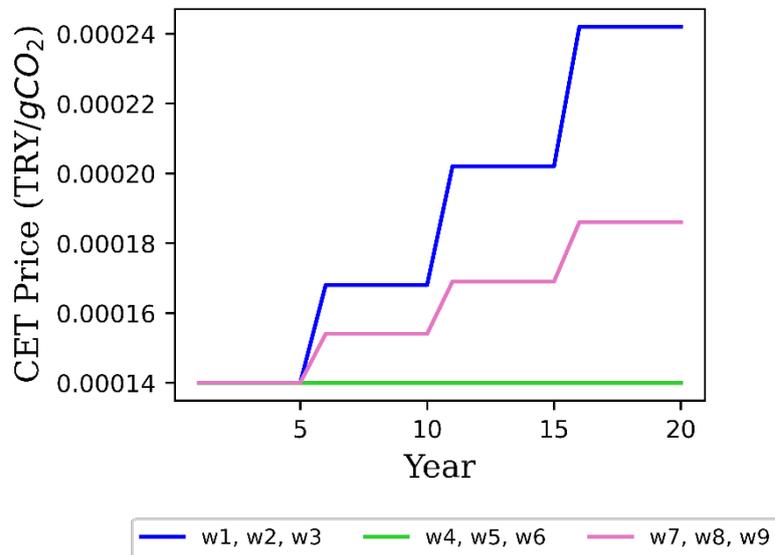

**Figure 11.** Carbon emission trading price for scenarios in Case 2.

The optimal equipment set for Case 2 is given in Table 7 with corresponding rated power and capacity values. The problem size is represented in Table 8. The drastic increase in the computational complicity due to the increase in the number of scenarios can be observed in Table 8 and Table 5.

Table 7. Chosen equipment for Case 2.

| Equipment | Rated Power (kW) | Capacity (kWh) |
|---|---|---|
| WindTur1 | 42961.02 | |
| PV1 | 177585.9 | |
| OilCogen1 | 8642.8 | |
| Battery | 23807.5 | 100000 |
| RecipEngine3 | 9341 | |
| FuelCell1 | 1400 | |
| BioGasif1 | 6600 | |

Table 8. Size and CPU time for Case 2.

| | Stochastic Solution |
|---|---|
| CPU Time (sec) | 1171526.8 |
| Continuous Variables | 1,304,928 |
| Discrete Variables | 103,698 |
| Number of Constraints | 2,052,919 |

Optimal solutions show that the electrolyzer-methanation reactor is not chosen to be installed. This shows that investing in this technology is still not optimal for the investigated location even if the uncertainties are considered. On the other hand, it should be noted that $H_2$ and $CO_2$ purchase is not allowed from outside of the grid. Accordingly, for the methanation reactor to operate, the electrolyzer must be operated to provide the required hydrogen. This would enable this compact form of equipment to operate together while utilizing the $CO_2$ produced within the energy hub. Accordingly, PtG is still one of the important technologies where the aim is to reduce the emission when $CO_2$ inside the grid is used. The sale of excess SNG term in the objective function would drive the model to produce and use SNG within the grid and sell the excess amount if it is profitable. In other words, if the income stream from the sale of excess SNG is less than the cost of installing and operating the PtG technology, then the installation decision of PtG technology would not be economically optimal. For that reason, there is no SNG production within the energy hub, and there is no income term from excess SNG trade. However, there is an income term from the sale of cap and trade, which is 30.9 million TRY corresponding to 0.3% of the total cost stream. The optimal cost for this case study is 11.5 billion TRY when a 12% discount rate is considered.

To observe the effect of the carbon emission trading price parameter, power balance figures from year 16 are given. Figs. 12-14 are the power balances for scenarios $w_1$, $w_4$, and $w_7$, respectively, belonging to year 16.

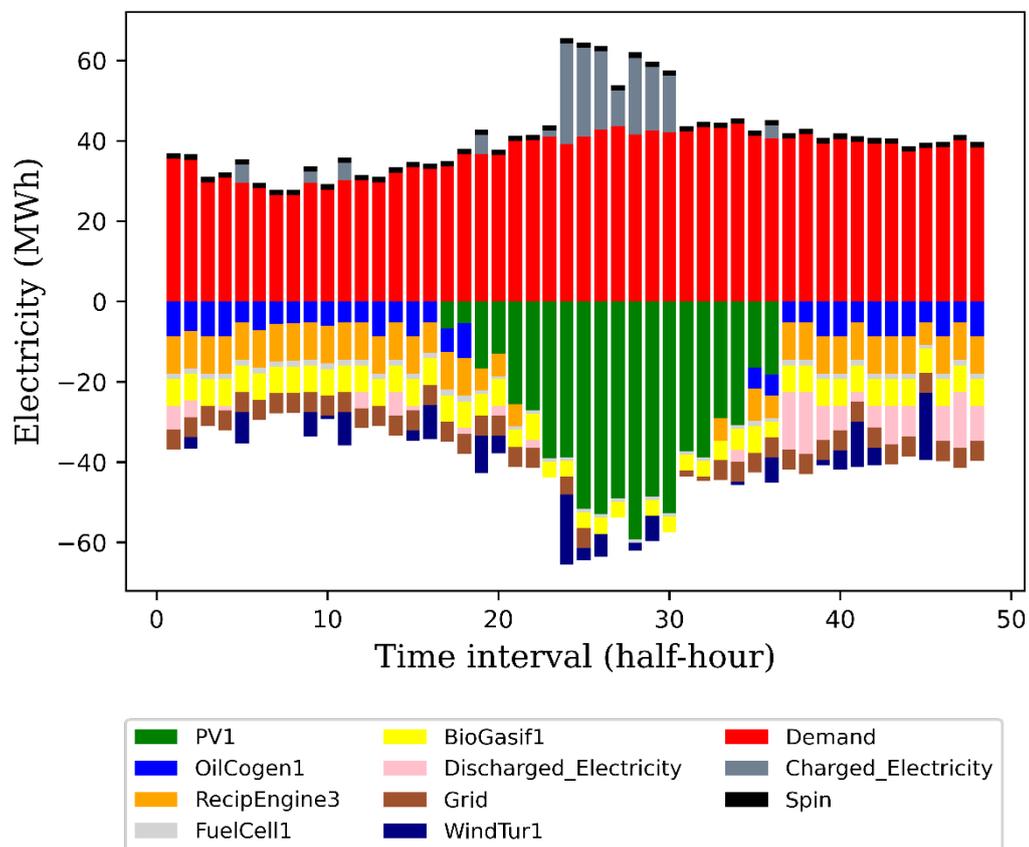

**Figure 12.** Optimal recourse results for a day in the 16$^{th}$ year for scenario $w_1$ in Case 2.

**Figure 13.** Optimal recourse results for a day in the 16$^{th}$ year for scenario $w_4$ in Case 2.

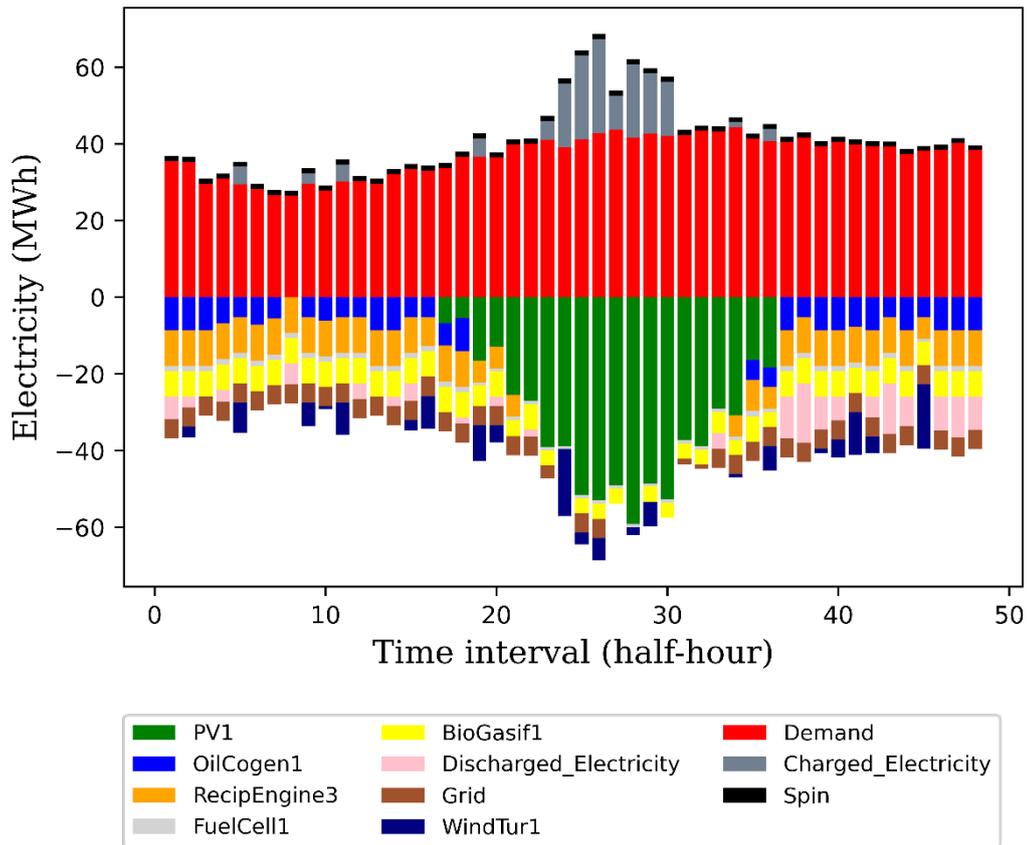

**Figure 14.** Optimal recourse results for a day in the 16$^{th}$ year for scenario $w_7$ in Case 2.

One thing in common in Figs. 12-14 is that most electricity is supplied with solar cells. It is followed by reciprocating combustion engines and biogasifiers, respectively. When we examine the total supplied electricity amount, it is observed that 53.6% of the power is supplied with renewable equipment. This percentage is found for scenarios $w_1$, $w_4$, and $w_7$, which include the likely case's demand data. In the same way, the power supply percentage from renewable energy-based equipment for scenario $w_1$ in Case 1, which has the likely case's demand data, is found as 52.9%. This is not a significant increase when Cases 1 and 2 are compared, but it shows the tendency towards renewable energy-based equipment selection.

### 4.3. Case 3

In Case 3, the $CO_2$ emission limit, *Elimit*, is determined as one of the uncertain parameters in addition to solar irradiance, wind speed, and ambient temperature, and the income term from the 'cap and trade' is not included in the objective function. The objective function only includes the income stream for the trade of excess SNG and the cost stream in the sub-case studies in Case 3 with the aim of minimizing net present cost.

Table 9. Uncertainties for Case 3.

| Scenarios | Electricity Demand | Operating Power Coefficient of Renewables as the Fraction of the Rated Power | | | $CO_2$ Emission Limit | Probability |
| --- | --- | --- | --- | --- | --- | --- |
| | | Wind Speed | Solar Irradiance | Temperature | | |
| $w_1$ | Likely Case | Likely Case | Likely Case | Likely Case | Decreasing | 1/9 |
| $w_2$ | Midlikely Case | Midlikely Case | Midlikely Case | Midlikely Case | Decreasing | 1/9 |
| $w_3$ | Unlikely Cases | Unlikely Case | Unlikely Case | Unlikely Case | Decreasing | 1/9 |
| $w_4$ | Likely Case | Likely Case | Likely Case | Likely Case | Constant | 1/9 |
| $w_5$ | Midlikely Case | Midlikely Case | Midlikely Case | Midlikely Case | Constant | 1/9 |
| $w_6$ | Unlikely Cases | Unlikely Case | Unlikely Case | Unlikely Case | Constant | 1/9 |
| $w_7$ | Likely Case | Likely Case | Likely Case | Likely Case | Slightly decreasing | 1/9 |
| $w_8$ | Midlikely Case | Midlikely Case | Midlikely Case | Midlikely Case | Slightly decreasing | 1/9 |
| $w_9$ | Unlikely Cases | Unlikely Case | Unlikely Case | Unlikely Case | Slightly decreasing | 1/9 |

It is mentioned that the $CO_2$ emission limit for a one-time interval (half an hour) is chosen as 280g$CO_2$/kWh as a starting point. The possible values for uncertain parameters that are used in scenario generation for case studies are given in Table 1. For this case study, it is predicted that the CO$_2$ emission limit will be constant for the first three scenarios. There will be a slight decrease in the CO$_2$ emission limit for the other three scenarios. Finally, for the remaining three scenarios, it is predicted that there will be a decrease in the CO$_2$ emission limit. These three cases are combined with the uncertain parameters of Case 1, and nine scenarios are created accordingly. It should be noted that the CET price is not included as an uncertain parameter in this case study. The decrease percentages will differ in Case 3.a and Case 3.b, which will be explained in the following subsections.

**Case 3.a**

The proposed scenarios regarding the CO$_2$ emission limit change are given in Fig. 15. For scenarios with a slight decrease in the CO$_2$ emission limit, there is a 10% drop every five years. On the other hand, for scenarios $w_1$, $w_2$ and $w_3$, there is a 20% decrease in the CO$_2$ emission limit every five years. These decrease percentages in CO$_2$ emission limit are chosen to be coherent with the regulations regarding Paris Agreement. For example, in scenarios with a 20% decrease in the CO$_2$ emission limit every five years, the drop from 280gCO$_2$/kWh to 143.36gCO$_2$/kWh in the limit is obtained at the end of the project lifetime, which corresponds to foreseeing approximately 50% drop in the CO$_2$ emission limit during the 20-year lifetime of the project.

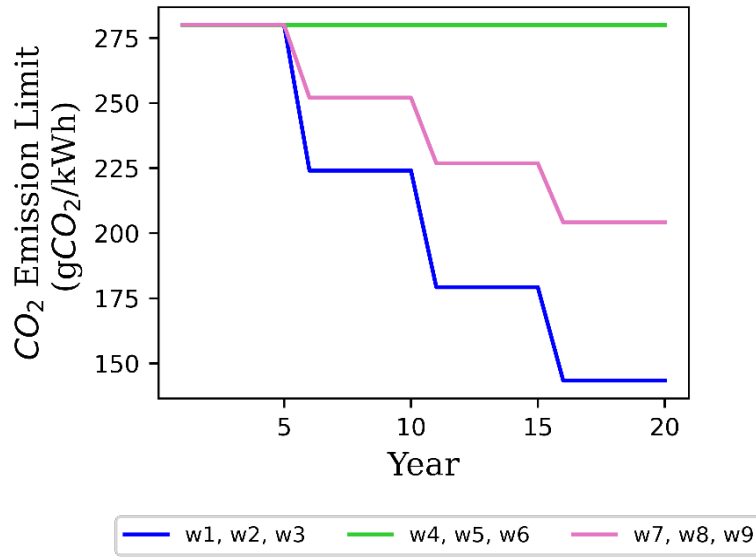

**Figure 15.** $CO_2$ emission limit for scenarios in Case 3.a.

The optimal results of Case 3.a suggest installing the equipment list given in Table 10 with corresponding rated power values. Additionally, the problem size and the CPU time for this case are given in Table 11. In Case 2, wind turbines, solar cells, and biogasifiers are chosen as renewable energy based equipment. On the other hand, in Case 3.a, a biomass generator is chosen additionally as a renewable energy based equipment that uses biomass.

**Table 10.** Chosen equipment for Case 3.a.

| Equipment | Rated Power (kW) | Capacity (kWh) |
|---|---|---|
| WindTur2 | 39560.13 | |
| PV2 | 189907.7 | |
| BioGen | 4349.02 | |
| OilCogen1 | 2343.8 | |
| Battery | 23533.2 | 100000 |
| RecipEngine3 | 9341 | |
| MicrTurb3 | 950 | |
| FuelCell1 | 1400 | |
| BioGasif1 | 6600 | |

Table 11. Size and CPU time for Case 3.a.

|  | Stochastic Solution |
|---|---|
| CPU Time (sec) | 154546.1 |
| Continuous Variables | 1,304,906 |
| Discrete Variables | 103,698 |
| Number of Constraints | 2,052,897 |

Since the methanation reactor and electrolyzer equipment pair are not chosen to be installed, there is no SNG production. Power balance figures for Case 3.a are given in Figs. 16-18 for scenarios $w_1$, $w_4$, and $w_7$ for year 16. These figures have the same electricity demand, solar irradiance, wind speed, and temperature data, which belong to Likely Case data. Hence, Figs. 16-18 have the same demand values and power supplied from wind turbines and solar cells, and they do not yield significantly different profiles.

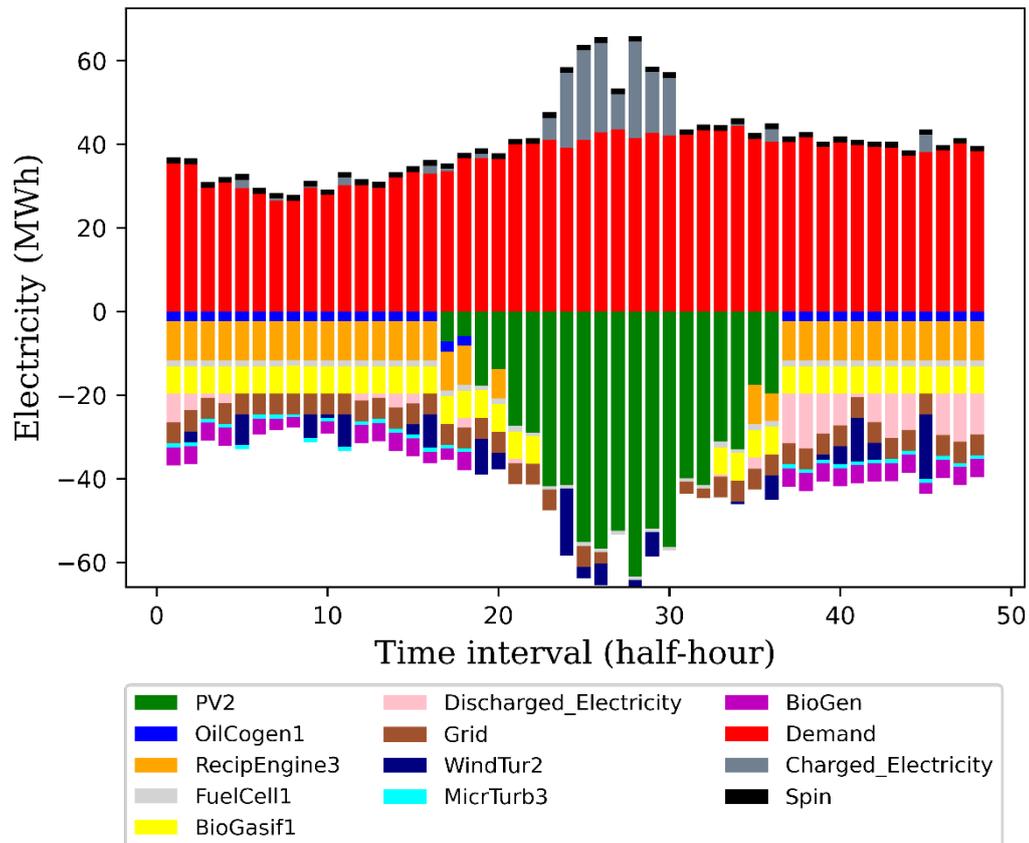

Figure 16. Optimal recourse results for a day in the 16$^{th}$ year for scenario $w_1$ in Case 3.a.

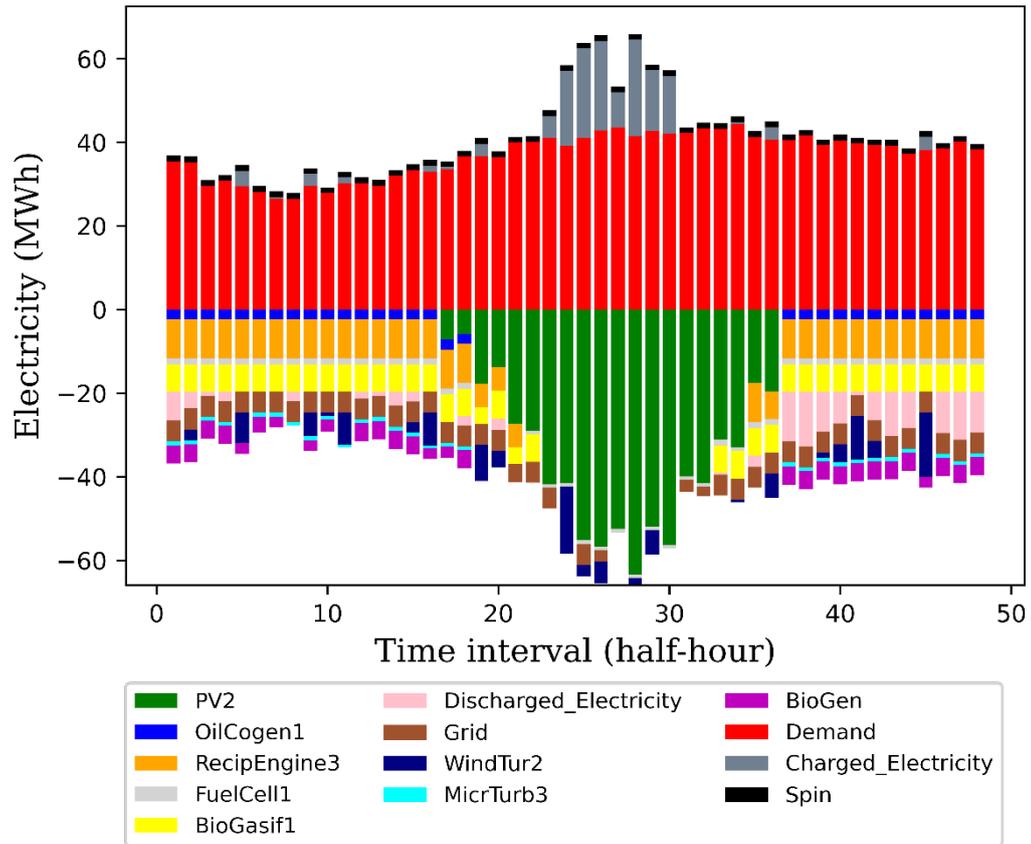

**Figure 17.** Optimal recourse results for a day in the 16th year for scenario $w_4$ in Case 3.a.

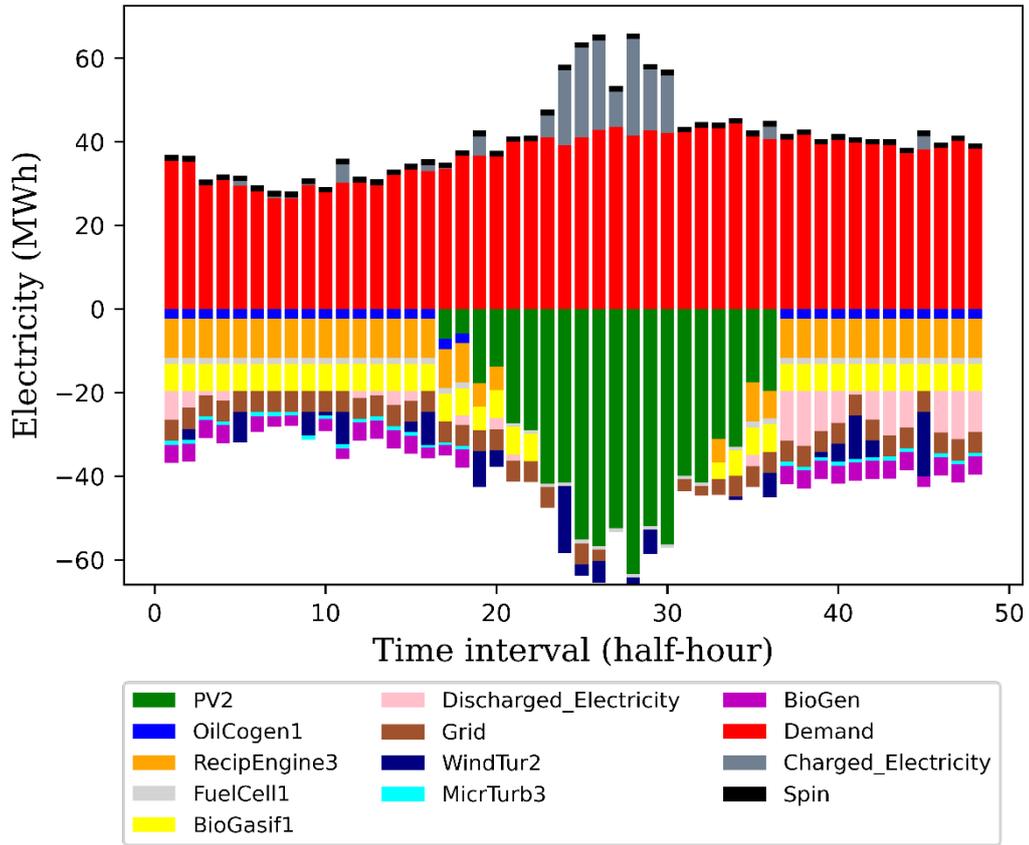

**Figure 18.** Optimal recourse results for a day in the 16th year for scenario $w_7$ in Case 3.a.

The renewable energy based equipment chosen to be installed for Case 3.a are wind turbines, solar cells, a biomass generator, and biogasifiers. In Fig. 19, the percentage of power supplied by renewable energy based equipment for every scenario for year 16 is given. For scenarios $w_1$, $w_4$, and $w_7$, the average percentage of renewable energy supply is 59.52 for Case 3.a. This percentage is 53.57 for Case 2. Thus, it can be yielded that, to adapt to the regulations from the Paris Agreement and global warming, renewable energy integration is one of the key components of the solution when several uncertainties are considered simultaneously for both optimal design and operation. Moreover, it is highly likely to expect a decrease in the $CO_2$ emission limit or an increase in CET price in the future government regulations that will try to prevent the worldwide $CO_2$ emission increase. At the same time, there will still be around 40% of the power source requirement from the non-renewable based generators. As a result, optimal integration and scheduling of the equipment must be studied under the impact of uncertainty and simultaneously.

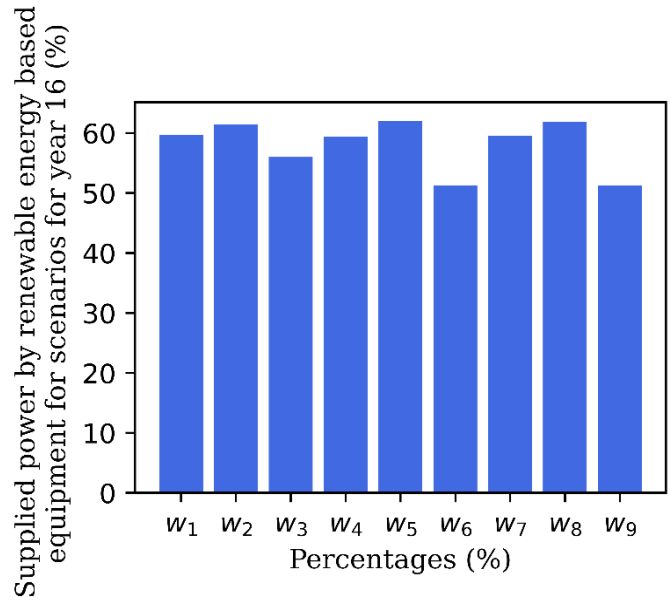

**Figure 19.** Power supply percentage from renewable energy-based equipment for scenarios in Case 3.a.

**Case 3.b**

In Case study 3.a, decrease percentages in $CO_2$ emission limit are defined as 10% for three scenarios and 20% for the other three scenarios, and finally, 0% for the last three scenarios every five years of the project lifetime, as stated previously. These percentages are chosen to be compatible with the emission regulations in Türkiye. Choosing steeper reduction rates than these decrease percentages would be considered unlikely at first glance for Türkiye. However, decrease percentages in $CO_2$ emission limit are chosen as 0%, 20%, and 50% for every five-year period of the project for Case 3.b to observe the effect of $CO_2$ emission limit if it is decreased further compared to Case 3.a. The 20% and 50 % decrease in $CO_2$ emission limit in every five-year period of the project corresponds to dropping the $CO_2$ emission limit to its half value for the three scenarios and nearly zero emission for the other three scenarios during the lifetime of the project.

The $CO_2$ emission limits for scenarios for Case 3.b are given in Figure 20. The optimal results obtained for Case 3.b are represented in Table 12.

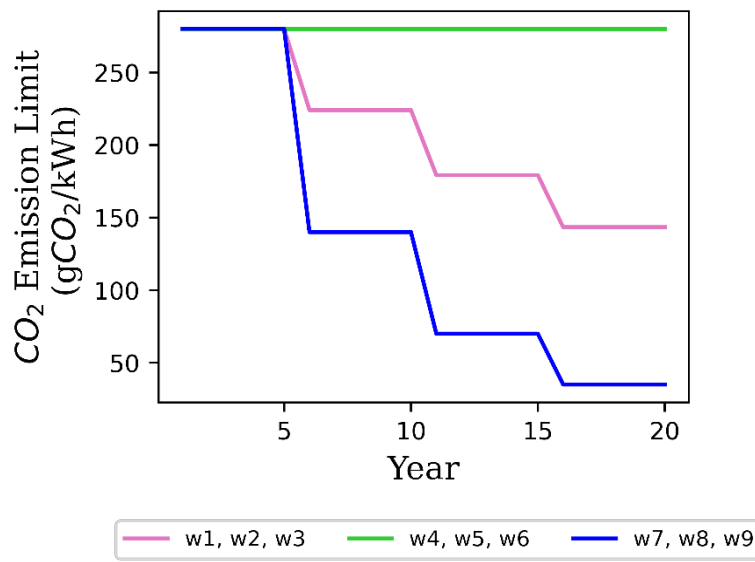

**Figure 20.** $CO_2$ emission limit for scenarios in Case 3.b.

**Table 12.** Chosen equipment for Case 3.b.

| Equipment | Rated Power (kW) | Capacity (kWh) |
|---|---|---|
| WindTur1 | 179494.3 | |
| PV2 | 181068.85 | |
| BioGen | 18558.45 | |
| Battery | 70298.45 | 137210.7 |
| RecipEngine3 | 9341 | |
| SteamTurb1 | 500 | |
| FuelCell1 | 1400 | |
| BioGasif1 | 6600 | |

When a higher decrease in $CO_2$ emission limit is considered, an increasing trend toward the selection of renewable energy based equipment is observed, as shown in Table 10 and Table 12. The results of Case 3.a in Table 10 suggest the installation of oil co-generator and micro turbines. In Case 3.b, instead of an oil co-generator and micro turbines, steam turbines are chosen to be installed. Additionally, it is observed that the rated power of the wind turbines given in Table 12 for this case study is way higher than the rated power of the same equipment in Case 3.a, as given in Table 10.

Although three of the scenarios foresee that there will be nearly zero emission limit in the 20-year period, the methanation reactor and electrolyzer equipment pair are not still optimal. It should be noted that the $CO_2$ emissions of wind turbines and solar cells are assumed to be zero. That is why when there is a higher decrease in $CO_2$ emission limit, the model prefers to install wind turbines with higher rated power values instead of choosing the PtG technology. Even though the price of the equipment of the PtG technology is enforced to be free, the equipment set is still not optimal. Hence, the reason for sub-optimality or infeasibility for green hydrogen and power-to-gas is not the installation cost of the equipment. PtG technology becomes cheaper when the $H_2$ is produced with an excess amount of renewable energy, such as wind turbines and solar cells, in a "green" way. Otherwise, operational costs enhance significantly. Therefore, convenient wind speed and solar irradiance profiles are required for the specific location in addition to available land.

To illustrate further, the optimality of PtG technology is analyzed using the deterministic model of this formulation. It is observed that the PtG technology is optimal when the following changes are made; the wind speed data is multiplied by five for every time interval, and the $CO_2$ emission limit is relaxed. In addition to these changes, the price of wind turbines, methanation reactor, and electrolyzer are enforced to be free, and the price of water as a raw material is made free of charge. After these changes are applied all together, the installation of PtG technology becomes optimal, which is evidently practically impossible.

A comparison of results of Case 3.a and Case 3.b in terms of net present costs of the project and the emitted $CO_2$ amount for a specific year are given in Table 14. The net present cost of the project for Case 3.b is increased by 35% approximately when compared to the net present cost of Case 3.a because of the significant decrease in $CO_2$ emission limits. When Case 3.a and 3.b are compared, in other words, a further decrease in $CO_2$ emission limit is considered for the scenarios, a drastic decrease in yearly emitted $CO_2$ is observed, which can be observed in Table 14.

**Case 3.c**

Finally, in this case study, PtG technology is forced to be operated by assigning the minimum operating power. The optimal equipment selection for this case study is given in Table 17. It is observed that solar cells, biomass generator, and biogasifiers are chosen to be installed as renewable energy based equipment together with the PtG technology. It is observed that the rated power of solar cells reaches nearly its maximum limit of 200 MW, which again shows the importance of green renewable energy for PtG technology.

The optimal objective function value for this case study is 28.1 billion TRY when a 12% discount rate is considered, and the electrolyzer is forced to operate at a minimum operating power of 20 MW for every time interval during the lifetime of the project.

Table 13. Chosen equipment for Case 3.c.

| Equipment | Rated Power (kW) | Capacity (kWh) |
|---|---|---|
| PV1 | 191699.1 | |
| PV2 | 200000 | |
| BioGen | 62582.2 | |
| Battery | 58659.64 | 100000 |
| RecipEngine3 | 9341 | |
| MicrTurb3 | 950 | |
| FuelCell1 | 1400 | |
| BioGasif1 | 6600 | |
| Electrolyzer | 45151.02 | |
| MetReac | 51300 | |

The net present cost values of the case studies and the yearly emitted $CO_2$ amount of the scenarios that belong to year 16 are given in Table 18. When other case studies in this chapter are compared, an increasing value of net present cost is obtained. This is directly related to the adapted $CO_2$ emission reduction strategies that result in the selection of more renewable energy-based equipment with higher rated power values.

To sum up, in case study 1, $CO_2$ emission is only limited according to the Paris Agreement. Additionally, in Case 2, carbon emission trading price is included as an uncertain parameter with an increasing trend in the scenarios. In Case 3, the $CO_2$ emission limit is obtained as an uncertain parameter. When the annual results in Table 18 are examined from Case 1 through Case 3.b, there is a decreasing trend in emission values of the corresponding scenarios. This is again directly related to the adopted $CO_2$ limiting strategies that limit the carbon emission further, starting from Case 1 to Case 3.b.

Table 14. Net present cost values and annual emitted $CO_2$ values of scenarios for year 16.

|  | Emission (tCO$_2$/year) | | | | |
| --- | --- | --- | --- | --- | --- |
| **Scenarios** | **Case 1** | **Case 2** | **Case 3.a** | **Case 3.b** | **Case 3.c** |
| $w_1$ | 144,447.3 | 124,936.7 | 90,333.8 | 48,771.4 | 46,128.7 |
| $w_2$ | 82,644.2 | 71,368.6 | 49,467.9 | 23,138.7 | 43,975.8 |
| $w_3$ | 134,259.6 | 122,670 | 80,017.8 | 57,986.6 | 54,127.7 |
| $w_4$ | - | 123,638.9 | 91,195.7 | 50,395.3 | 45,645.4 |
| $w_5$ | - | 70,873.2 | 49,550.6 | 23,231.4 | 44,034.2 |
| $w_6$ | - | 122,166.7 | 93,986.4 | 58,036.9 | 53,229.9 |
| $w_7$ | - | 124,827.3 | 91,157.4 | 19,649.6 | 22,647.8 |
| $w_8$ | - | 71,161.9 | 49,162.4 | 16,844.5 | 17,737.2 |
| $w_9$ | - | 123,994.6 | 93,986.2 | 19,282.8 | 19,535.6 |
| **Net present Cost (billion TRY)** | 11.1 | 11.5 | 11.8 | 15.9 | 28.1 |

## 5. Conclusion

This work proposes a multi-period two-stage stochastic programming model in which six uncertain parameters are considered, and the different combinations of these are investigated in various case studies. Three base cases are used, namely likely, mid-likely, and unlikely. The effect of uncertain demand, wind speed, temperature, and solar irradiance data are investigated. Electricity demand, wind speed, temperature, solar irradiance data, carbon emission trading (CET) price, and $CO_2$ emission limit are considered uncertain parameters. The model is formulated to meet one-third of the power and heat demand of Yalova city while minimizing the net present cost for 20 years of plant life, where the cost stream is compensated with two income streams: the sale of 'cap and trade' and excess SNG. An increasing and decreasing trend in the possible values of CET price and $CO_2$ emission limit is adapted to follow the regulations regarding Paris Agreement.

Cases 2 and 3 include uncertain parameters of CET price and $CO_2$ emission limit in addition to the uncertain parameters in Case 1. When the first stage installation decision variables are compared for Case 1, Case 2, and Case 3.a, a trend toward selecting renewable energy-based equipment with higher rated power values is observed. The average percentage of power supplied with renewable energy based equipment is found as 52.9%, 53.6%, and 59.52% for Case 1, Case 2, and Case 3.a, respectively, when the scenarios that use the likely case's demand data are compared. Although the percentages for Cases 1 and 2 are quite close to each other, it still states that the model tends to choose more renewable energy based equipment when CET price with the increasing profile is included as an uncertain parameter. While adapting greater decrease in $CO_2$ emission limit in Case 3.b results in 15.9 billion TRY when compared with Case 3.a, which results in 11. 8 billion TRY corresponding to 35% increase approximately. Enforcing the installation of PtG technology in Case 3.c results in 28.1 billion TRY net present cost, which corresponds to 77% increase when it is compared with the net present cost of Case 3.b. When the $CO_2$ emission values of the scenarios are investigated separately, starting from Case 1 to Case 3.c, a decrease in annual $CO_2$ emission is observed while obtaining higher net present costs. This is because there is a trend of installing more renewable energy-based equipment with increasing rated power values. That is the reason for the reduced emission trend and rising net present costs from Case 1 to Case 3.c because renewable energy-based equipment is relatively expensive.

While Case 2 considers the CET price as an additional uncertain parameter, Case 3 considers the $CO_2$ emission limit as an additional uncertain parameter. When the annual emitted $CO_2$ rates of scenarios of Case 2 and 3 are compared, it is observed that the annual carbon emission of Case 3 is lower than Case 2. This can be interpreted as implementing the $CO_2$ emission limit as an uncertain parameter instead of the CET price would be more beneficial in terms of fighting against global warming. Of course, the assigned $CO_2$ limits in scenarios play an important role in this, but the strict limits should be chosen anyway for the scenarios.

Additionally, the candidate pool of Cases 2 and 3 includes PtG technology in addition to the candidate pool of Case 1, where the optimality of the power to gas technology is investigated in the second and third cases. The PtG technology is not optimal for this study mainly because of the available renewable energy resources for this specific location in Türkiye and since the PtG becomes optimal when $H_2$ production is cheaper with excess wind and solar energy.

**Acknowledgments**

This publication has been produced benefiting from the 2232 International Fellowship for Outstanding Researchers Program of TUBITAK (Project No: 118C245). However, the entire responsibility of the publication belongs to the owner of the publication.**REFERENCES**

[1] BBC News. Temperature hits 34.8C on Scotland's hottest day n.d. https://www.bbc.com/news/uk-scotland-62225963 (accessed July 28, 2022).

[2] CNN World. The UK's hottest day destroyed their homes. They fear it's a sign of worse ahead n.d. https://edition.cnn.com/2022/07/26/uk/london-fire-dagenham-heat-wave-climate-gbr-cmd-intl/index.html (accessed July 28, 2022).

[3] BBC Türkçe. Almanya enerji krizine çare arıyor: "Daha kısa duş alın" n.d. https://www.bbc.com/turkce/articles/ckvx941883wo?at_custom1=%5Bpost+type%5D&at_custom4=7A50E126-09AD-11ED-8383-FEF7923C408C&at_campaign=64&at_custom3=BBC+Turkey&at_custom2=facebook_page&at_medium=custom7&fbclid=IwAR0VBUEup42NFIsyI0EgxFZcFRegerTBGT-ufD2VwQ (accessed July 28, 2022).

[4] BBC News. Climate change: EU unveils plan to end reliance on Russian gas 2022. https://www.bbc.com/news/science-environment-60664799#:~:text=The plan envisages ending reliance,of renewables%2C biogas and hydrogen. (accessed May 4, 2022).

[5] Davis SJ, Lewis NS, Shaner M, Aggarwal S, Arent D, Azevedo IL, et al. Net-zero emissions energy systems. Science (80- ) 2018;360. https://doi.org/10.1126/science.aas9793.

[6] De Luca G, Fabozzi S, Massarotti N, Vanoli L. A renewable energy system for a nearly zero greenhouse city: Case study of a small city in southern Italy. Energy 2018;143:347–62. https://doi.org/10.1016/j.energy.2017.07.004.

[7] Li X, Damartzis T, Stadler Z, Moret S, Meier B, Friedl M, et al. Decarbonization in Complex Energy Systems: A Study on the Feasibility of Carbon Neutrality for Switzerland in 2050. Front Energy Res 2020;8:1–17. https://doi.org/10.3389/fenrg.2020.549615.

[8] Farrokhifar M, Nie Y, Pozo D. Energy systems planning: A survey on models for integrated power and natural gas networks coordination. Appl Energy 2020;262:114567. https://doi.org/10.1016/j.apenergy.2020.114567.

[9] Dolatabadi A, Mohammadi-Ivatloo B, Abapour M, Tohidi S. Optimal Stochastic Design of Wind Integrated Energy Hub. IEEE Trans Ind Informatics 2017;13:2379–88. https://doi.org/10.1109/TII.2017.2664101.

[10] Zhou Z, Zhang J, Liu P, Li Z, Georgiadis MC, Pistikopoulos EN. A two-stage stochastic programming model for the optimal design of distributed energy systems. Appl Energy 2013;103:135–44. https://doi.org/10.1016/j.apenergy.2012.09.019.

[11] Yu J, Ryu JH, Lee I beum. A stochastic optimization approach to the design and operation planning of a hybrid renewable energy system. Appl Energy 2019;247:212–20. https://doi.org/10.1016/j.apenergy.2019.03.207.

[12] Iris Ç, Lam JSL. Optimal energy management and operations planning in seaports with41

# SUPPLEMENTARY MATERIAL

## A.1. Mathematical Formulation of the Model

**Nomenclature**

    **Definition of sets;**

| | |
|---|---|
| $I$ | : Set of equipment |
| $G$ | : Set of generators ($G \in I$) |
| $S$ | : Set of storage units ($S \in I$) |
| $R$ | : Set of renewables ($R \in I$) |
| $K$ | : Set of years |
| $T$ | : Set of time-intervals |
| $N$ | : Set of resources |
| $W$ | : Set of uncertainty (Discrete) |
| $w$ | : Scenarios |

    **Definition of first-stage decision variables;**

| | |
|---|---|
| $rp_i$ | : Rated power of equipment $i$, in kW |
| $b_i$ | : Capacity of storage unit $i$, in kWh ($i \in S$) |
| $a_i$ | : Binary decision to install an equipment or not |

    **Definition of second-stage decision variables;**

| | |
|---|---|
| $sp_{nkw}$ | : Spinning reserve of a generator in scenario $w$ ($i \in G$) |
| $kc_{iktw}$ | : Binary decision to run a generator or storage at interval $l$ in a day of year $k$ in scenario $w$ ($i \in G \cup S$) |
| $p_{iktw}$ | : Operating power of generator $i$ at interval $l$ in a day of year $k$ in scenario $w$, in kW ($i \in G$) |
| $pch_{iktw}$ | : Charging power of equipment $i$ at interval $l$ in a day of year $k$ in scenario $w$, in kW ($i \in S$) |
| $pdch_{iktw}$ | : Discharging power of equipment $i$ at interval $l$ in a day of year $k$ in scenario $w$, in kW ($i \in S$) |
| $soc_{ik0w}$ | : Initially stored energy of a storage device |

  in a day of year *k* in scenario *w*, in kWh ($i \in S$)

$soc_{iktw}$  : Stored energy of a storage device

  in a day of year *k* in scenario *w*, in kWh ($i \in S$)

$yx_{nktw}$  : Surplus output of resource *n* at time interval *t* in a day of year *k* in scenario *w*, (kWh-MJ-g)/hour

$u_{nktw}$  : System input of resource *n* at time interval *t*

  in a day of year *k* in scenario *w*, (kWh-MJ-g)/hour

$\xi_{iktw}$  : Peak penalty

**Definition of Parameters;**

$int$  : interest rate

$inf$  : inflation rate

$D$  : 365 days per year

$\Delta T$  : Time interval, in hour

$N_I$  : Maximum number of introducible equipment

$LCO2_{kw}$  : Daily $CO_2$ emission limit, $gCO_2$/kWh

$rp_i^{min}$  : Minimum introducible rated power of equipment *i*, in kW

$rp_i^{max}$  : Maximum introducible rated power of equipment *i*, in kW

$b_i^{min}$  : Minimum introducible capacity of storage *i*, in kWh ($i \in S$)

$b_i^{max}$  : Maximum introducible capacity of storage *i*, in kWh ($i \in S$)

$p_{iktw}^{min}$  : Minimum operating power of equipment *i* at time interval *t* in a day of year *k* by the fraction of the rated power, dimensionless ($i \in G$)

$p_{iktw}^{max}$  : Maximum operating power of equipment *i* at time interval *t* in a day of year *k* by the fraction of the rated power, dimensionless ($i \in G$)

$q_{iktw}^{min}$  : Minimum operating charged energy of a storage *i* at time interval *t* in a day of year *k* by the fraction of the rated power, dimensionless ($i \in S$)

$q_{iktw}^{max}$  : Maximum operating charged energy of a storage *i* at time interval *t* in a day of year *k* by the fraction of the rated power, dimensionless ($i \in S$)

$g_i$  : Generation of resource *n* when equipment *i* is run at the unit power (kW) for unit time (hour)

$c_i$  : Consumption of resource *n*

  when equipment *i* is run at the unit power (kW) for unit time (hour)

$d_{nktw}$  : Demand of resource *n* at time step *t* in a day of year *k*

$\bar{p}_{iktw}$  : Operating power coefficient of the renewable equipment *i*

  at time step *t* in a day of year *k* ($i \in R$)

| | |
|---|---|
| $\alpha_i^0$ | : Cost for building equipment *i* per unit power, |
| $\beta_i^0$ | : Cost for building equipment *i* per unit capacity ($i \in S$) |
| $\gamma_i^0$ | : Fixed cost for building equipment *i* |
| $\alpha_i^k$ | : Cost for maintenance of equipment *i* per unit power in year *k* |
| $\beta_i^k$ | : Cost for maintenance of equipment *i* per unit capacity in year *k* ($i \in S$) |
| $\gamma_i^k$ | : Fixed cost for maintenance of equipment *i* in year *k* |
| $\phi_{nktw}^+$ | : Cost for unit system output of resource *n* at time interval *t* in a day of year *k* |
| $\phi_{nktw}^-$ | : Cost for unit system input of resource *n* at time interval *t* in a day of year *k* |
| $P_w$ | Probability of occurrence of scenario *w* |
| $\delta$ | : peak penalty coefficient |
| $Elimit_{kw}$ | $CO_2$ emission limit of scenario w, in year k |
| $Eprice_{kw}$ | Carbon emission trading (CET) price of scenario *w*, in year *k* |
| $\eta$ | : efficiency of the photovoltaic cell |
| $\Phi$ | : solar irradiance, in kW |
| $\kappa$ | : temperature correction factor for the photovoltaic cell |
| $T_p$ | : ambient temperature, in $^0C$ |
| $T_{p,ref}$ | : reference temperature, in $^0C$ |
| $v_{cut,in}$ | : cut in wind speed, in m/sec |
| $v_{cut,out}$ | : cut out wind speed, in m/sec |
| $v_{rated}$ | : rated wind speed, in m/sec |
| $v$ | : actual wind speed, in m/sec |

Objective function:

$$\min \quad f^{Initial}(rp, b, a) + \sum_{k \in \mathcal{K}} f_k^{Operational}(rp, b, a, yx_{knlw}, u_{knlw}) - \sum_{k \in \mathcal{K}} f_k^{Cap\ and\ Trade} \quad (1)$$

$$- \sum_{k \in \mathcal{K}} f_k^{Excess\ SNG}$$

Eq. 1 represents the objective function which is the summation of initial investment cost and operating & maintenance costs in addition to two income streams over the lifetime of the project. The objective function aims to minimize the net present cost. The detailed representation of these cost terms is given in Eqs. 17 and 18. The cost stream is compensated by two income terms. In order words, two income terms have a reducing effect on the overall cost. These income streams come from selling excess SNG (synthetic natural gas) and carbon allowances, referred to as cap and trade.

Installation and sizing decisions:

$$\sum_{i \in I} a_i \leq N_I, \quad (2)$$

$$rp_i^{min} a_i \leq rp_i \leq rp_i^{max} a_i, \quad i \in G, \quad (3)$$
$$b_i^{min} a_i \leq b_i \leq b_i^{max} a_i, \quad i \in S, \quad (4)$$
$$a_i \in \{0,1\}, \quad i \in I.$$

The total number of equipment that can be installed is limited by ten because of the area restrictions, which is shown in Eq. 2. *a* term in Eqs. 3 and 4 refer to binary decision variable to install a piece of equipment or not. The *rp* and *b* terms in Eqs. 3 and 4 refer to the rated power of generators and the capacity values for batteries, respectively. These equations represent the rated power and capacity constraints when the decision of installation is made, which corresponds to an *a* value of 1.

Unit Commitment and Operational Constraints:

$$p_{iktw}^{min}(rp_i - (1 - kc_{iktw})rp_i^{max}) \leq p_{iktw} \leq p_{iktw}^{max} rp_i, \quad (5)$$

$$i \in G, k \in \mathcal{K}, t \in \mathcal{T}, w \in W$$

$$0 \leq p_{iktw} \leq p_{iktw}^{max} rp_i^{max} kc_{iktw} \quad , \quad i \in G, k \in \mathcal{K}, t \in \mathcal{T}, w \in W \tag{6}$$

$$kc_{iktw} \in \{0,1\}, \quad i \in G, k \in \mathcal{K}, t \in \mathcal{T}, w \in W$$

While the *rp* term corresponds to the rated power of the equipment, the *p* term corresponds to the operating power of the generators. Eqs. 5 and 6 represent the upper and lower limit of the operating power values of the generators, where the $kc_{iktw}$ term represents the binary decision to operate a generator or not, which corresponds to switched on and off mode. It should be noted that except for the first stage variables, second stage variables have a value for every scenario *w*.

$$0 \leq pch_{iktw} \leq p_{iktw}^{max} rp_i \quad i \in S, k \in \mathcal{K}, t \in \mathcal{T}, w \in W \tag{7}$$

$$0 \leq pdch_{iktw} \leq p_{iktw}^{max} rp_i \quad i \in S, k \in \mathcal{K}, t \in \mathcal{T}, w \in W \tag{8}$$

$$pch_{iktw} \leq p_{iktw}^{max} rp_i^{max} ks_{iktw} \quad i \in S, k \in \mathcal{K}, t \in \mathcal{T}, w \in W \tag{9}$$

$$pdch_{iktw} \leq p_{iktw}^{max} rp_i^{max}(1 - ks_{iktw}) \quad i \in S, k \in \mathcal{K}, t \in \mathcal{T}, w \in W \tag{10}$$

$$ks_{iktw} \in \{0,1\} \quad i \in S, k \in \mathcal{K}, t \in \mathcal{T}, w \in W$$

Eq. 7 represents the charging power limit, where $pch_{iktw}$ is the charging power of storage equipment, and Eq. 8 represents discharging limit, where $pdch_{iktw}$ is discharging power of storage equipment. It should be noted that $ks_{iktw}$ represents the binary decision to operate storage equipment or not, and $p_{iktw}^{max}$ corresponds to the maximum operating power by the fraction of the rated power for every piece of equipment. Eqs. 9 and 10 represent that the storage equipment cannot be charged and discharged simultaneously.

Storage Constraints:

$$q_{iktw}^{min} b_i \leq soc_{iktw} \leq q_{iktw}^{max} b_i \quad i \in S, k \in \mathcal{K}, t \in \mathcal{T}, w \in W \tag{11}$$

$$soc_{ik|T-1|w} = soc_{ik0w} \quad i \in S, k \in \mathcal{K}, w \in W \tag{12}$$

$$soc_{iktw} = \begin{cases} soc_{ik0w} + \Delta T(pch_{iktw} - pdch_{iktw}) & \text{if } t = t(1) \\ soc_{ikt-1w} + \Delta T(pch_{iktw} - pdch_{iktw}) & \text{otherwise} \end{cases} \tag{13}$$

$$i \in S, k \in \mathcal{K}, t \in \mathcal{T}, w \in W$$

$soc_{iktw}$ term in Eqs. 11-13 represents the reserved energy in storage equipment. $q_{iktw}^{min}$ and $q_{iktw}^{max}$ terms correspond to the minimum and maximum operating charged energy of storage equipment which are 0.2 and 0.8, respectively. Hence, Eq. 11 states that the reserved energy in storage equipment should be between 20% and 80% of the storage capacity. Eq. 12 represents that the storage equipment's initial and final conditions must be the same on a particular day. Also, Eq. 13 represents the calculation of the $soc_{iktw}$ term with the charging and discharging of storage equipment.

Non-negativity Constraints for System Inputs and Outputs:

$$u_{nktw} \geq 0 \quad n \in \mathcal{N}, k \in \mathcal{K}, t \in \mathcal{T}, w \in W \tag{14}$$

Eq. 14 represents the non-negativity constraint of the system input of resources, $u_{nktw}$.

$$yx_{nktw} \geq 0 \quad n \in \mathcal{N}, k \in \mathcal{K}, t \in \mathcal{T}, w \in W \tag{15}$$

Eq. 15 represents the non-negativity constraint of the surplus output of resource n, $yx_{nktw}$.

Power Balance:

$$\sum_{i \in G} g_{in} p_{iktw} + \sum_{i \in S} g_{in} pdch_{iktw} + \sum_{i \in \mathcal{R}} g_{in} \bar{p}_{iktw} rp_i + u_{nktw} = \sum_{i \in G} c_{in} p_{iktw}$$
$$+ \sum_{i \in S} c_{in} pch_{iktw} + \sum_{i \in \mathcal{R}} c_{in} \bar{p}_{iktw} rp_i + sp_{nkw} + yx_{nktw} + d_{nktw} \tag{16}$$

$$n \in \mathcal{N}, k \in \mathcal{K}, t \in \mathcal{T}, w \in W$$

Eq. 16 represents the material balance for every resource *n*. It includes the generation, consumption, system input, surplus output, demand, and spin of a resource *n*. Also, $g_{in}$ and $c_{in}$ terms correspond to the generation and consumption of resources $n$ by the equipment. The values for $g_{in}$ and $c_{in}$ are given in Appendix Tables A.3 and A.4. The terms $sp_{nkw}$ and $d_{nktw}$ correspond to spinning reserve and demand, respectively. The $sp_{nkw}$ must be less than three percent of the maximum demand during the day. It should be noted that the $g_{in}$ and $c_{in}$ in Appendix are taken from the literature. These tables represent the amount of generated or consumed resource $n$ by equipment $i$ operating at the unit power (kW) for unit time (hour).

Cost Calculations:

$$f^{initial}(rp, b, a) = \sum_{i \in E}(\alpha_i^0 rp_i + \gamma_i^0 a_i) + \sum_{i \in S} \beta_i^0 b_i, \quad (17)$$

$$f_k^{Operational} = (\sum_{i \in E}(1 + inf)^{k-1}(\alpha_i^k rp_i + \gamma_i^k a_i)$$

$$+ \sum_{i \in S}(1 + inf)^{k-1} \beta_i^k b_i)$$

$$+ D\Delta T \sum_{t \in T} \sum_{n \in \mathcal{N}} \sum_{w \in W}(1 + inf)^{k-1} P_w (\phi_{nktw}^+ yx_{nktw} + \phi_{nktw}^- u_{nktw})$$

$$+ \sum_{w \in W} P_w \, D|\mathcal{K}||T|\xi \quad (18)$$

$$f_k^{Cap \, and \, Trade} = D\Delta T \sum_{w \in W} P_w \, Eprice_{kw} \, (LCO2_{kw} - E_{actual})$$

$$= D\Delta T \sum_{w \in W} P_w \, Eprice_{kw} \, (LCO2_{kw} - \sum_{i \in T}(yx_{,CO2\,,k,t,w})) \quad (19)$$

$$f_k^{Excess \, SNG} = D\Delta T \sum_{t \in T} \sum_{w \in W}(1 + inf)^{k-1} P_w \left(0.90 * \phi_{,Gas,ktw}^- yx_{,Gas,k,t,w}\right) \quad (20)$$

$$\xi \geq \delta(p_{gktw})(g_{,Electricity,kt}) \quad (21)$$

Eqs. 17 and 18 represent the cost streams of the model, including the initial cost and total operational and maintenance costs (Q&M). The rated power and capacity values of the equipment are multiplied with the corresponding cost term to find the initial and total maintenance costs, where inflation is also considered. It should be noted that producing excess resource, $n$, is represented with the $\phi^+_{nktw} yx_{nktw}$ term, which is neglected in this study. Additionally, the $\phi^-_{nktw}$ term represents the cost of system input of resource, $n$, resulting in the $\phi^-_{nktw} u_{nktw}$ term corresponding to the raw material cost. The last term in Eq. 18 represents the penalty cost associated with peak shaving. The $D$ and $|T|$ terms correspond to the number of days in a year and the size of time intervals, respectively.

It should be noted that the $CO_2$ emission is limited according to Paris Agreement for this study. Additionally, an income stream is added to the objective function coming from the sale of carbon allowances, namely 'cap and trade,' represented in Eq. 19. The term multiplied with *Eprice* is the gap between carbon allowance and actual carbon emission. This gap is multiplied by carbon emission trading (CET) price, *Eprice*, and probability of occurrence resulting in the income stream coming from the cap and trade. It should be noted that the *Eprice* corresponds to the previously mentioned uncertain parameter of carbon emission trading (CET) price, which is introduced as an uncertain parameter only in Case 2.

The second income term comes from selling excess synthetic natural gas (SNG) produced in the energy hub represented in Eq. 20. Excess SNG is multiplied by the cost of natural gas together with the probability of occurrence. It should be noted that the natural gas price is taken 90% of the purchased gas price for the energy grid while also considering the inflation.

Eq. 21 represents the peak penalty selection equation. While the $\delta$ term corresponds to the peak penalty coefficient, the $\xi$ term corresponds to the peak penalty. Eq. 21 states that the peak penalty must be greater than the maximum electricity amount produced from the generators. This provides shaving the peak supply from the generators. It should be noted that smaller values are chosen for the peak penalty coefficients for this study.

Emission Constraints:

$$\sum_{i \in T}(yx_{,CO2\,,k,t,w}) \leq LCO2 \qquad (22)$$

$$LCO2 = Elimit * \sum_{i \in T} Power\ Load_{ktw} \qquad (23)$$

It will be mentioned in the paper later that the $CO_2$ emission from renewable energy based equipment is assumed as zero. Hence, when Eq.16 is examined for the $CO_2$, it is observed that the $yx_{,CO2\,,k,t,w}$ term corresponds to generated $CO_2$ minus the consumed $CO_2$ emission inside the energy hub. Including the $CO_2$ tax in the objective function and paying penalties for excess production of $CO_2$ without limiting the $CO_2$ emission would not serve the aforementioned challenges in the energy industry in terms of economic feasibility and global warming concerns. Therefore, the net $CO_2$ emission is limited daily in which the corresponding limit for one (half an hour) time interval is 280gCO$_2$/kWh. It means that the $CO_2$ emission could be higher or less than this limit for the one-time interval, but the daily limit is less than (48)(280) gCO$_2$/kWh.

Eq. 22 states that the net emitted $CO_2$ amount is limited. This limit corresponds to the emission limit of a one-time interval multiplied by the power load. Hence, the $LCO2$ term corresponds to the carbon emission allowance, in other words, the carbon emission quote of the energy grid [24]. It should be noted that the *Elimit* in Eq. 23 corresponds to the previously mentioned uncertain parameter of $CO_2$ emission limit, which is introduced as an uncertain parameter only in Case 3.

Limiting $CO_2$ emissions is basically stating that the emission should be less than a certain amount for the model. Adding the sale of cap and trade as an income stream in objective function together with emission limitations only allows trading the cap; it does not allow excess $CO_2$ emission. At the same time, trading the carbon allowance with this emission limitation also reduces the overall cost with the income stream from selling the cap. Hence, adding the term of trading $CO_2$ cap to the objective function would drive the model to emit less $CO_2$ further with the aim of minimizing the net present cost. In other words, emitting less $CO_2$ results in a greater amount of carbon allowance that can be traded, resulting in a greater income stream that further reduces the net present cost. Hence, it can be stated that this situation restricts the model further to follow the current emission regulations and help to minimize the cost further simultaneously.

Renewable Equipment Generations:

$$\bar{p}_{\prime PV\prime ktw} = \eta \Phi \left(1 - \kappa(T_{p,w} - T_{p,ref})\right) \tag{24}$$

$$\bar{p}_{\prime WT\prime ktw} = \begin{cases} 0 & if\ v < v_{cut,in} \\ 0 & if\ v > v_{cut,out} \\ 1 & if\ v_{rated} \leq v \leq v_{cut,out} \\ \dfrac{v - v_{cut,in}}{v_{rated} - v_{cut,in}} & if\ v_{cut,in} \leq v \leq v_{rated} \end{cases} \tag{25}$$

Eq. 24 calculates the solar cells' operating power coefficient, which is proportional to the efficiency of the solar cell, $\eta$, and the solar irradiance, $\Phi$. In the same way, Eq. 25 represents the operating power coefficient of wind turbines. $v$ terms correspond to wind speed which should be between $v_{cut,in}$ and $v_{cut,out}$ values for the wind turbine to operate. Also, the $v_{rated}$ term corresponds to the wind speed value in which the wind turbines operate at maximum power [9].

It should be noted that the operating power coefficient of solar cells and wind turbines are given as the operating power of renewables as the fraction of rated power of the renewable equipment. In other words, the multiplication of the operating power coefficient and rated power of renewables gives the operating power of renewables. The operating power coefficient of renewables given in Eqs. 24 and 25 is one of the uncertain parameters that depends on the uncertain parameters of wind speed, solar irradiance, and temperature for three cases. The corresponding values of these uncertain parameters are given in Section 3, in Figs. 1-5.

## A. Tables

**Table A.1.** Power and Capacity Limits [26–29]

|  | Minimum Rated Power [kW] | Maximum Rated Power [kW] | Minimum Capacity | Maximum Capacity |
|---|---|---|---|---|
| Wind Turbine-1 | 100 | 200000 | 0 | 0 |
| Wind Turbine-2 | 100 | 200000 | 0 | 0 |
| Wind Turbine-3 | 100 | 200000 | 0 | 0 |
| Photovoltaic-1(PV) | 100 | 200000 | 0 | 0 |
| Photovoltaic-2(PV) | 100 | 200000 | 0 | 0 |
| Photovoltaic-3(PV) | 100 | 200000 | 0 | 0 |
| Biomass Generator | 2500 | 500000 | 0 | 0 |
| Gas Co-generator-1 | 3350 | 670000 | 0 | 0 |
| Oil Co-generator -1 | 750 | 150000 | 0 | 0 |
| Lithium Ion Battery | 1000 | 500000 | 100000 kWh | 500000 kWh |
| Reciprocating Internal Combustion Engine-3 | 9341 | 9341 | 0 | 0 |
| Gas/Combustion Turbine-1 | 3304 | 3304 | 0 | 0 |
| Steam Turbine-1 | 500 | 500 | 0 | 0 |
| Micro Turbine-3 | 950 | 950 | 0 | 0 |
| Fuel Cell-1 | 1400 | 1400 | 0 | 0 |
| Biogasifier-1 | 6600 | 6600 | 0 | 0 |
| Electrolyzer | 10000 | 100000 | 0 | 0 |
| Methanation Reactor | 51300 | 51300 | 0 | 0 |

**Table A.2.** Initial and Maintenance Costs For The First Year [26–29]

|  | Initial Cost | | Maintenance Cost | |
|---|---|---|---|---|
| Symbol | Per Rated Power [TL/kW] | Per Capacity [TL/kWh-MJ] | Per Rated Power [TL/kW] | Per Capacity [TL/kWh-MJ] |
| Wind Turbine-1 | 20039.04 | 0 | 403.2 | 0 |
| Wind Turbine-2 | 20039.04 | 0 | 403.2 | 0 |
| Wind Turbine-3 | 20039.04 | 0 | 403.2 | 0 |
| Photovoltaic-1(PV) | 20744.64 | 0 | 248.6 | 0 |
| Photovoltaic-2(PV) | 20744.64 | 0 | 248.6 | 0 |
| Photovoltaic-3(PV) | 20744.64 | 0 | 248.6 | 0 |
| Biomass Generator | 28082.88 | 0 | 1814.4 | 0 |
| Gas Co-generator-1 | 8537.76 | 0 | 672 | 0 |
| Oil Co-generator -1 | 9172.8 | 0 | 530.9 | 0 |
| Lithium Ion Battery | 0 | 12096 | 134.4 | 336 |
| Reciprocating Internal Combustion Engine-3 | 10031 | 0 | 514.1 | 0 |
| Gas/Combustion Turbine-1 | 22967 | 0 | 762.05 | 0 |
| Steam Turbine-1 | 7952 | 0 | 604.8 | 0 |
| Micro Turbine-3 | 17500 | 0 | 725.8 | 0 |
| Fuel Cell-1 | 32200 | 0 | 2419.2 | 0 |
| Biogasifier-1 | 34392 | 0 | 2693.4 | 0 |
| Electrolyzer | 5355 | 0 | 289 | 0 |
| Methanation Reactor | 5580.1 | 0 | 210.6 | 0 |

**Table A.3.** Consumption Table at Unit Power [kW] In Unit Time [h].[26–29]

| | Consumption | | | | | | | | | |
|---|---|---|---|---|---|---|---|---|---|---|
| Equipment Name | Elect. [kWh/kW/h] | Heat [MJ/kW/h] | Biomass [MJ/kW/h] | Gas [kWh/kW/h] | Oil [MJ/kW/h] | $CO_2$ [$gCO_2$/kW/h] | Wood F. [MJ/kW/h] | Coal [MJ/kW/h] | $H_2$ [kWh/kW/h] | Water [kg/kW/h] |
| Wind Turbine-1 | 0 | 0 | 0 | 0 | 0 | 0 | 0 | 0 | 0 | 0 |
| Wind Turbine-2 | 0 | 0 | 0 | 0 | 0 | 0 | 0 | 0 | 0 | 0 |
| Wind Turbine-3 | 0 | 0 | 0 | 0 | 0 | 0 | 0 | 0 | 0 | 0 |
| Photovoltaic-1(PV) | 0 | 0 | 0 | 0 | 0 | 0 | 0 | 0 | 0 | 0 |
| Photovoltaic-2(PV) | 0 | 0 | 0 | 0 | 0 | 0 | 0 | 0 | 0 | 0 |
| Photovoltaic-3(PV) | 0 | 0 | 0 | 0 | 0 | 0 | 0 | 0 | 0 | 0 |
| Biomass Generator | 0 | 0 | 20.07 | 0 | 0 | 0 | 0 | 0 | 0 | 0 |
| Gas Co-generator-1 | 0 | 0 | 0 | 3.81 | 0 | 0 | 0 | 0 | 0 | 0 |
| Oil Co-generator -1 | 0 | 0 | 0 | 0 | 11.3 | 0 | 0 | 0 | 0 | 0 |
| Lithium Ion Battery | 1.05 | 0 | 0 | 0 | 0 | 0 | 0 | 0 | 0 | 0 |
| Reciprocating Internal Combustion Engine-3 | 0 | 0 | 0 | 2.4 | 0 | 0 | 0 | 0 | 0 | 0 |
| Gas/Combustion Turbine-1 | 0 | 0 | 0 | 4.18 | 0 | 0 | 0 | 0 | 0 | 0 |
| Steam Turbine-1 | 0 | 0 | 0 | 0 | 57.4 | 0 | 0 | 0 | 0 | 0 |
| Micro Turbine-3 | 0 | 0 | 0 | 3.76 | 0 | 0 | 0 | 0 | 0 | 0 |
| Fuel Cell-1 | 0 | 0 | 0 | 2.35 | 0 | 0 | 0 | 0 | 0 | 0 |
| Biogasifier-1 | 0 | 0 | 0 | 0 | 0 | 0 | 20.5 | 0 | 0 | 0 |
| Electrolyzer | 1.32 | 0 | 0 | 0 | 0 | 0 | 0 | 0 | 0 | 0 |
| Methanation Reactor | 0 | 0 | 0 | 0 | 0 | 177.14 | 0 | 0 | 1.28 | 0.748 |

**Table A.4.** Generation Table at Unit Power [kW] In Unit Time [h]. [26–29]

| Equipment Name | Generation | | | | | | | | | |
|---|---|---|---|---|---|---|---|---|---|---|
| | Elect. [kWh/kW/h] | Heat [MJ/kW/h] | Biomass [MJ/kW/h] | Gas [kWh/kW/h] | Oil [MJ/kW/h] | $CO_2$ [$gCO_2$/kW/h] | Wood F. [MJ/kW/h] | Coal [MJ/kW/h] | $H_2$ [kWh/kW/h] | Water [kg/kW/h] |
| Wind Turbine-1 | 1 | 0 | 0 | 0 | 0 | 0 | 0 | 0 | 0 | 0 |
| Wind Turbine-2 | 1 | 0 | 0 | 0 | 0 | 0 | 0 | 0 | 0 | 0 |
| Wind Turbine-3 | 1 | 0 | 0 | 0 | 0 | 0 | 0 | 0 | 0 | 0 |
| Photovoltaic-1(PV) | 1 | 0 | 0 | 0 | 0 | 0 | 0 | 0 | 0 | 0 |
| Photovoltaic-2(PV) | 1 | 0 | 0 | 0 | 0 | 0 | 0 | 0 | 0 | 0 |
| Photovoltaic-3(PV) | 1 | 0 | 0 | 0 | 0 | 0 | 0 | 0 | 0 | 0 |
| Biomass Generator | 1 | 0 | 0 | 0 | 0 | 79 | 0 | 0 | 0 | 0 |
| Gas Co-generator-1 | 1 | 5.7 | 0 | 0 | 0 | 599 | 0 | 0 | 0 | 0 |
| Oil Co-generator -1 | 1 | 2.5 | 0 | 0 | 0 | 738 | 0 | 0 | 0 | 0 |
| Lithium Ion Battery | 0.95 | 0 | 0 | 0 | 0 | 0 | 0 | 0 | 0 | 0 |
| Reciprocating Internal Combustion Engine-3 | 1 | 3.02 | 0 | 0 | 0 | 448.2 | 0 | 0 | 0 | 0 |
| Gas/Combustion Turbine-1 | 1 | 6.32 | 0 | 0 | 0 | 361.5 | 0 | 0 | 0 | 0 |
| Steam Turbine-1 | 1 | 41.9 | 0 | 0 | 0 | 247.6 | 0 | 0 | 0 | 0 |
| Micro Turbine-3 | 1 | 4.93 | 0 | 0 | 0 | 369.7 | 0 | 0 | 0 | 0 |
| Fuel Cell-1 | 1 | 3.33 | 0 | 0 | 0 | 235.9 | 0 | 0 | 0 | 0 |
| Biogasifier-1 | 1 | 0 | 0 | 0 | 0 | 106.5 | 0 | 0 | 0 | 0 |
| Electrolyzer | 0 | 0 | 0 | 0 | 0 | 0 | 0 | 0 | 1 | 0 |
| Methanation Reactor | 0.002 | 0 | 0 | 1 | 0 | 0 | 0 | 0 | 0 | 0 |